\documentclass[]{compositio}

\usepackage{amsbsy,amssymb,amscd,amsfonts,latexsym,amstext,delarray,amsmath,stackrel,lineno,tabularx,url,MnSymbol,tikz,cite,yfonts}
\usepackage{float}
\usepackage[pdfborder={0 0 0}]{hyperref}
\usepackage{mathrsfs}
\usepackage[all]{xy}
\hypersetup{
    colorlinks = true,
    linkcolor = blue,
    anchorcolor = red,
    citecolor = blue,
    filecolor = red,
    pagecolor = red,
    urlcolor = magenta
}

\usepackage{graphicx}

\usepackage{pdfpages}

\usepackage{blindtext}
\usepackage{etoolbox}
\newcommand{\printthis}[2][true]{%
\ifbool{#1}{%
#2%
}{%
}%
}

\newtheorem{thm}{Theorem}[section] 
\newtheorem{defn}[thm]{Definition} 
\newtheorem{prop}[thm]{Proposition}

\newtheorem{lem}[thm]{Lemma}

\newtheorem{rem}[thm]{Remark}

\def\fii{{\rm  finite}}

\def\Aut{{\rm Aut}}

\def\End{{\rm End}}

\def\Hom{{\rm Hom}}
\def\id{{\rm id}}

\def\SL{{\rm SL}}

\def\Spec{{\rm Spec\,}}

\def\tr{{\rm tr}}

\def\Mo{\mathfrak{Mo}}
\def\A{{\mathbb A}}

\def\C{{\mathbb C}}
\def\F{{\mathbb F}}
\def\K{{\mathbb K}}
\def\N{{\mathbb N}}

\def\Q{{\mathbb Q}}
\def\R{{\mathbb R}}
\def\Z{{\mathbb Z}}

\def\W{{\mathbb W}}

\def\tr{{\rm tr}}

\def\cA{{\mathcal A}}

\def\cC{{\mathcal C}}

\def\cE{{\mathcal E}}

\def\cH{{\mathcal H}}

\def\cF{{\mathcal F}}

\def\cM{{\mathcal M}}
\def\cO{{\mathcal O}}
\def\cP{{\mathcal P}}

\def\cR{{\mathcal R}}
\def\cS{{\mathcal S}}

\def\qqq{\,,\,~\forall}

\def\pik{\pi^{\rm new}}

\newcommand{\ie}{{\it i.e.\/}\ }

\newcommand{\opcit}{{\it op.cit.\/}\ }

\def\W{{\mathbb W}}
\def\gh{{\rm gh}}

\def\id{{\mbox{Id}}}

\def\dim{{\mbox{dim}}}

\def\Hom {{\mbox{Hom}}}

\def\End{{\mbox{End}}}

\def\sss{{\mathbb S}}
\def\gop{{\Gamma^{\rm op}}}

\def\zmax{{\Z_{\rm max}}}

\def\fr{{\rm Fr}}

\def\sspec{{\pmb{\mathfrak{Spec}}}}

\def\Se{\frak{ Sets}}
\def\fin{\frak{ Fin}}

\def\An{\mathfrak{ Ring}}
\def\spzb{{\overline{\Spec\Z}}}

\def\vvert{{\Vert}}

\def\sdd{{\rm Sd}^*}
\def\dop{{\Delta^{\rm op}}}

\def\lbt{{\tilde \Lambda}}

\def\epi{{(\tilde\Lambda^{\rm op})^\wedge}}

\def\nt{\N^{\times}}
\def\wnt{{\widehat{\N^{\times}}}}

\def\Ses{{\Se_*}}
\def\crel{{\Se_{2,*}}}
\def\sses{{\cS_*}}

\def\bfh{{\bf H}}

\def\zmm{\hat\Z^\times}
\def\Mo{\mathfrak{Mo}}
\def\An{\mathfrak{ Ring}}
\def\Ab{\mathfrak{ Ab}}
\def\Mr{\mathfrak{ MR}}

\def\salg{\mathfrak{ S}}

\def\Se{\mathfrak{ Sets}}

\def\smod{{\sss-{\rm Mod}}}

\begin{document}

\title{ Segal's Gamma rings and\\ universal arithmetic}
\author{Alain Connes}
\email{alain@connes.org}
\address{College de France, I.H.E.S. and Ohio State University}
\author{Caterina Consani}
\email{kc@math.jhu.edu}
\address{Department of Mathematics, The Johns Hopkins
University\newline Baltimore, MD 21218 USA}
\dedication{In memoriam Michael Atiyah,  with profound admiration}
\thanks{We are grateful to Jens Hemelaer for his useful comments on an earlier version of this paper.}

\begin{abstract}

 Segal's $\Gamma$-rings provide a natural framework for absolute algebraic geometry. We use  Almkvist's global Witt construction to explore the relation with J. Borger $\F_1$-geometry and compute the  Witt functor-ring $\W_0(\sss)$ of the simplest $\Gamma$-ring $\sss$.  We  prove that it is isomorphic to the Galois invariant part of the BC-system, and exhibit the close relation between $\lambda$-rings and the Arithmetic site. Then, we  concentrate on the Arakelov compactification $\spzb$ which  acquires a structure sheaf of $\sss$-algebras. After supplying a probabilistic interpretation of the classical theta invariant of a divisor $D$ on $\spzb$, we show how to associate to $D$ a $\Gamma$-space  that encodes, in homotopical terms, the Riemann-Roch problem for $D$.
\end{abstract}

\maketitle

\vspace{0.1in}

\section{Introduction}

Michael Atiyah was quite right stating that ``without dreams there is
no art, no mathematics no life''. Dreaming new paths in mathematics
has helped mathematicians to overcome seemingly impossible technical obstacles,  frequently  guiding them to discover  fascinating
new theories. It has brought us to explore the enchanted realm of J. Tits's  ``field with 1 element'' \cite{Tits}
and to develop in \cite{schemeF1} the theory of the ``absolute point" Spec($\F_1$) and the algebraic
geometry over it (speculated in \cite{Manin}). This is the topic pursued in this article that we heartily dedicate to
M.~F.~Atiyah in honor of his
farsighted intuition and ability to communicate enthusiasm  while spreading new ideas.\newline
In his seminal paper \cite{Tits}, Tits, motivated by his geometric interpretation  of  C. Chevalley's work on algebraic groups of Lie type, advocated the search of a mysterious ``field with one element'' $\F_1$ (denoted $K_1$ in that paper), that would also account for the understanding of  the infinite series of simple alternating groups $\text{Alt}(n)$ ($n\geq 5$), as $\SL_n(\F_1)$ \cite{KS}.    In the last thirty years, this search was pursued in a number of ways: we refer, in particular, to \cite{KS,Soule,Borgers,TV,Durov}. From a quite different perspective, the algebraic $K$-theory of spaces initiated by F. Waldhausen led to  the notion of ``brave new rings"  thus shifting the interest of algebraic topologists from rings to ring spectra, where the integers $\Z$  become an algebra over the sphere spectrum \cite{Goodwillie}. A suitable model for connective spectra was provided by G. Segal's $\Gamma$-spaces which, endowed with Lydakis' smash product, form a closed, symmetric  and monoidal category.  This category is the natural framework for both Hochschild and cyclic homologies \cite{DGM}.   In  our research, we found that  the theory of {\em discrete} Segal's $\Gamma$-rings  provides a natural extension of the classical theory of rings and  semirings,  hence supplying a convincing algebraic incarnation for the $\F_1$-dream of a characteristic-free (universal) theory. To be more specific, a Segal's discrete  $\Gamma$-ring (we abbreviate it as $\Gamma$-ring for simplicity) is by definition a monoidal object in the category $\smod$ of pointed covariant functors from  finite pointed sets to pointed sets. In order to work with a small category one uses the skeleton category  $\gop$ with one object $n_+=\{0,\ldots, n\}$ ($0$ is the base point) for each integer $n$.  $\smod$ is a  closed, symmetric and monoidal category,  thus the monoids in it give a very concrete notion of  $\Gamma$-rings,  forming the category $\salg$.  The simplest $\Gamma$-ring $\sss$ corresponds to the identity functor and it is the most basic incarnation of the sphere spectrum over which the ring of integers  becomes an algebra (we shall use freely the terminology ``$\sss$-algebras" in place of $\Gamma$-rings).\newline
 Many attempts to model Tits' alleged ``field with one element" find naturally their place  among $\Gamma$-rings \cite{CCprel}. This holds for N. Durov's theory  \cite{Durov} and for the constructions with hyperrings, while both the category $\Mo$ of commutative pointed monoids (advocated in \cite{KOW, Kato, deit, TV}), and the category of semirings (that we used extensively in \cite{CCas,CCscal1,CCscal2,CCscal3,CCscal4} to model Tits' idea of ``characteristic one") are full subcategories of  $\salg$.
The connection with the ideas  developed by  Kapranov  and Smirnov in \cite{KS} is provided by the fact that $\smod$  can be seen  as the category of vector spaces over $\F_1$ as in \opcit, provided one works in the presheaf topos $\Gamma^\wedge$ on Segal's category $\Gamma$, rather than in the topos of sets \cite{CCgromov}. \newline
The relation with the original constructions of C. Soul\' e \cite{Soule} came as a pleasant surprise. Our earlier developments  \cite{compositio}  promoted, after the initial approach in \cite{ak}, a theory   of schemes (of finite type) over $\F_1$  using the category $\Mr$ obtained by gluing together the category $\Mo$ with the category  $\An$ of commutative rings. In  \cite{schemeF1} we realized that the category $\Mr$ is simply the full subcategory of $\salg$ whose objects are either in $\Mo$ or $\An$!
Thus, the a priori artificial process of glueing $\Mr=\An\cup_{\beta,\beta^*} \Mo$ described in \cite{compositio} is now fully justified and has suggested to us, in view of the results of \opcit  the development of algebraic geometry working directly in $\salg$.\newline  This is the path we have started to follow in \cite{schemeF1}, as far as the affine case is concerned, by extending to the general case of an arbitrary $\sss$-algebra $A$,  the construction of its spectrum as a {\em topos} $\sspec(A)$, endowed with a  structure sheaf of $\sss$-algebras. The spectrum $\sspec(A)$  derives from  a {\em Grothendieck site} endowed with a {\em presheaf} of $\sss$-algebras: this datum is more refined than that of the associated topos and sheaf.   When $A=HR$ for a semiring $R$, the associated topos-spectrum determines  the prime spectrum, while it hands   Deitmar's spectrum \cite{deit} when $A=\sss M$ is the spherical algebra associated to a monoid $M$. It is by now clear that our first development  \cite{compositio} is a special case of this new ``absolute algebraic geometry", and in particular that
$\sspec(\sss)$ is the natural candidate for $\Spec(\F_1)$.\newline
The distinctive feature of this  new theory is to give significance  to operations which are meaningless in ordinary algebraic geometry,  such as taking the quotient of a ring by a subgroup of its multiplicative group, or the restriction to the unit ball in  normed rings. These operations  make sense in full generality for $\sss$-algebras, thus one can analyze their effect on the associated spectra. \newline
The connection between this development and the theory of T\"oen-Vaqui\' e \cite{TV} (a general theory of algebraic geometry  for any symmetric monoidal closed category that is complete and cocomplete, like the category $\smod$) shows that the theory in \cite{TV}, when implemented to the category $\smod$, does {\em not} agree with ordinary algebraic geometry, already in the simplest case of the two point space corresponding to the spectrum of the product of two fields.

In the present paper we examine the connection of  \cite{schemeF1} with J. Borger's  $\F_1$-theory \cite{Borgers}. In Section \ref{sectwitt}, we compute the global Witt ring, as understood by Almkvist in \cite{Al,Al1}, of $\sss$ and we show that it is the $\lambda$-ring obtained from the integral BC-system by taking the fixed points of the Galois action of $\zmm$
\begin{thm}\label{w0sintro} The  global Witt ring $\W_0(\sss)$ is canonically isomorphic to the invariant part
$\Z[\Q/\Z]^G$ of the group ring $\Z[\Q/\Z]$, under the action of the  group $G=\Aut(\Q/\Z)= \widehat \Z^\times$.
\end{thm}
This result is in agreement with the idea developed in \cite{ArithBC}  that the integral BC-system is the Witt ring of the algebraic closure of $\F_1$. In fact we conjecture that the above Witt construction applied to  $\sss[\mu_\infty]$, where  $\mu_\infty$ denotes the group of abstract roots of unity with zero added, is canonically isomorphic to the full BC-system.
The key idea advocated in \cite{Borgers} is that $\lambda$-rings are obtained by extension of scalars from $\F_1$ to $\Z$. The basic result which greatly simplifies the definition of $\lambda$-rings is due to  C. Wilkerson \cite{Wilkerson} and formulates, in the flat case, the structure in terms of the Adams operations and Frobenius lifts, in place of the $\lambda$-operations themselves, whose algebraic rules are more complicated to state. There is a rather obvious relation between this formulation of $\lambda$-rings and the Arithmetic Site \cite{CCas} that we discuss in \S\ref{sectlambda}. Namely, a  $\lambda$-ring automatically generates, using the Adams operations, a sheaf of rings over the Arithmetic topos $\wnt$. In Section \ref{sectwitt}, we also recall the origin of  $\wnt$ from the cyclic and epicyclic categories and we discuss the interrelations of the various toposes candidates for the ``absolute geometric point".   \newline
Global fields  are either function fields in finite characteristic or number fields. The analogy between these two families breaks down at the
archimedean places for the lack of an ``absolute arithmetic" allowing to treat them on the same footing as the ultrametric places. The simplest case to consider is the field $\Q$  of the rational numbers. At a non-archimedean place $p$, the subset $\{x\in \Q_p\mid \vert x\vert \leq 1\}$ is a subring of the local field $\Q_p$ but the  corresponding set, at the archimedean place, does not define a subring of the local field $\R$. It does however still define  an $\sss$-subalgebra of $H\R$.  In \cite{CCprel} we applied the   theory of $\sss$-algebras to extend the structure sheaf of  $\Spec \Z$ to the Arakelov compactification $\overline{ \Spec \Z}$, as a sheaf of $\sss$-algebras.
The remaining sections of this paper  concentrate on the Arakelov compactification $\overline{ \Spec \Z}$  as an arithmetic curve over $\sss$ (in fact more precisely over $\sss[\mu_2]$). Each Arakelov divisor provides a natural sheaf of modules over this extended structure sheaf of $\overline{ \Spec \Z}$. Moreover, this new structure for $\overline{\Spec\Z}$ over $\sss$  endorses also a one parameter group of weakly invertible sheaves, whose tensor product rules are the same as   the composition rules   of  Frobenius correspondences over
the Arithmetic Site \cite{CCarith, CCas}. \newline
For a  function field  with field of constants $\F_q$, the dimension $\dim_{\F_q}H^0(D)$ of the $\F_q$-vector space $H^0(D)$ of solutions of the Riemann-Roch problem for a divisor $D$  is the logarithm $\log_q(\# H^0(D))$ of the  cardinality of $H^0(D)$.
Following the  analogy between function fields and number fields, $\dim_{\F_q}H^0(D)$ is replaced, in the number field case, by  the theta dimension $h_\theta^0(D)$ of an Arakelov divisor $D$ \cite{vdG-S, bost}. In Section \ref{secttheta} we introduce, in view of the new understanding of $\overline{\Spec\Z}$ over $\sss$, the probabilistic interpretation of $h_\theta^0(D)$, as the logarithm of an {\em average value}  of the  naive counting of solutions formulated in terms of  $\sss$-modules. More precisely, for a given  rank one discrete subgroup $L$ in a one dimensional {\em complex}  vector space $E$, one measures the length of vectors $\xi\in E$ by  comparing $\xi$ with the unit of length provided by $L$. This process yields the  non negative integer
  \begin{equation}\label{countingLzintro}
 [\xi/L]:=\# \{\ell \in L\mid \Vert \ell \Vert \leq \Vert \xi \Vert \}
\end{equation}
that is independent of the choice of the hermitian norm on $E$. Then, in Theorem~\ref{thmcountingLz} we show the equality
 \begin{equation}\label{averagecintro}
\exp(h_\theta^0(D))=\int [\xi/L] dp(\xi)
\end{equation}
where the  probability measure $dp(\xi)$ is  a Gaussian and depends on the choice of the  hermitian norm on $E$ through the normalization condition  $\int \Vert \xi \Vert dp(\xi)=\frac 12$.
  Our probabilistic interpretation of $h_\theta^0(D)$ is implicit in some proofs of inequalities involving the naive counting.  On the other hand, it differs substantially from the related work of J.~B.~Bost \cite{bost} based on the thermodynamic expression of the logarithm of a sum as a supremum involving the entropy function.  The role of the Gaussian on the complex vector space $E$ also suggests the existence of a possible  relation with the formulas of V.~Mathai and D.~Quillen  \cite{MQ} for the Chern character.\newline
  In  \cite{CCgromov} we found that the  $\Gamma$-ring arising as the stalk of the structure sheaf of $\overline{ \Spec \Z}$ at the archimedean place, is  intimately related to hyperbolic geometry and  Gromov's norm. This development involves a generalization of the homology of a simplicial complex when the coefficients are no longer abelian groups but  $\sss$-modules. This new homology takes values in $\sss$-modules and strictly extends the standard theory.  Here, the key fact is that the category of $\Gamma$-spaces, that  is the central tool of  \cite{DGM} and forms  the core of the local structure of algebraic $K$-theory, is the category of simplicial objects in $\smod$.
  Thus, according to the general principle inherent to the
Dold-Kan correspondence, homotopy theory of $\Gamma$-spaces ought to be considered as   absolute homological algebra of $\sss$-modules.   Since the category $\smod$ is not abelian,  the tools of homological algebra need to be replaced  along the line of the Dold-Kan correspondence,  that provides,  for an abelian category $\cA$  the correspondence between chain complexes in non-negative degrees and simplicial objects \ie objects of the category $\cA^\dop$.   While in \cite{DGM} the systematic use of fibrant replacements frequently simplifies constructions at the price of neglecting important information,  for our needs it is necessary   to work directly with $\sss$-modules rather than abelian groups, as the most genuine receptacle of homotopy invariants. This was the strategy already adopted in \cite{CCgromov}. \newline
Therefore  the  natural goal to pursue, in parallel with the elaboration of ``absolute algebraic geometry'',  is the construction of  the  homological algebra  over $\salg$.  The clear advantage of working with $\Gamma$-rings is  that this category forms the natural groundwork where cyclic homology is rooted \cite{DGM} and  in particular, where de Rham theory is naturally available.   \newline
 In Section \ref{sectdks}  we construct a first instance of the analogue of the Dold-Kan correspondence  for $\sss$-modules.
     We consider the adelic proof of the Riemann-Roch formula  as developed by A. Weil for global fields of positive characteristic. The analogue of his construction for the field $\Q$ provides one with a short complex $C(D)$ of $\sss$-modules  associated to an Arakelov divisor $D$ on $\overline{ \Spec \Z}$. This complex  arises from a natural morphism  of $\sss$-modules  canonically  associated to  $D$ (we refer to \S\ref{sectmorphi} for its definition). In Proposition \ref{propgam}, we construct the $\Gamma$-space $\bfh(D)$ associated to the complex $C(D)$ and hence to the Arakelov divisor $D$.  We view $\bfh(D)$  as  the homological incarnation of the Riemann-Roch problem for  $D$ and  in \S\S\ref{sectlevel1} and  \ref{secthigher} we  determine its homotopy.
 It is an interesting open question to find the analogue of Serre's duality in this context,  connecting (a version of) $\pi_0(\bfh(D))$ with $\pi_1(\bfh(-D))$. Obviously, the goal is to prove the Riemann-Roch formula $h_\theta^0(D)-h_\theta^0(-D)=\deg(D)$ for the probabilistic theta dimension in purely geometric terms,  without appealing to the analytic Poisson's formula.


\section{The global Witt ring $\W_0(\sss)$  }\label{sectwitt}

We recall that  $\sss$-modules (equivalently $\Gamma$-sets) are by definition covariant functors $\gop\longrightarrow\Ses$ between pointed categories.  Here, $\gop$ is   the small, full subcategory   of  the category $\fin_*$ of  finite pointed sets, whose objects are  the pointed sets $k_+:=\{0,\ldots ,k\}$, for each integer $k\geq 0$ ($0$ is the base point) and the morphisms connecting two objects are the sets $\gop(k_+,m_+)=\{f: \{0,1,\ldots,k\}\to\{0,1,\ldots,m\}\mid f(0)=0\}$.
The morphisms  in the category $\smod$ of  $\sss$-modules are natural transformations. The category $\smod$ is a closed, symmetric and monoidal category. \newline
An $\sss$-algebra is a monoid in $(\smod,\wedge,\sss)$:  the natural inclusion $\sss: \gop \longrightarrow \Ses$ is the simplest example. \newline
In this section we compute Almkvist's invariant $\W_0$ \cite{Al,Al1} for the simplest $\sss$-algebra, \ie we determine  $\W_0(\sss)$, motivated by the guess in \cite{Borgers} that $W(\F_1)=\Z$.\newline

\subsection{ Almkvist's invariant  $\W_0(\sss)$}\label{sectw0}
 In Almkvist's original definition of  the invariant $\W_0(A)=\tilde K_0(\text{end}\,{\bf P}(A))$  for a commutative ring $A$, one uses exact sequences to derive the notion of  additive invariants  for endomorphisms of finitely generated projective $A$-modules in the Grothendieck group $K_0(\text{end}\,{\bf P}(A))$. In particular, the condition stating that the  invariants are additive on exact sequences suffices to eliminate the non-trivial  part of the Jordan decomposition of a matrix,  thus it plays an essential role  in the  construction. By definition,  $\W_0(A)$  is the quotient of $K_0(\text{end}\,{\bf P}(A))$ by the subgroup $K_0(A)$ of trivial endomorphisms. In the context of  the non abelian category of $\sss$-modules  a generalization of Almkvist's construction is still possible. It suffices to require the additivity for exact sequences of endomorphisms  of $\sss$-modules of the form
\begin{equation}\label{exse}
	*\to T(E)\to E\to E/T(E)\to *\qquad \forall T\in \End(E),
\end{equation}
where the $\sss$-module $E/T(E)$ is defined by collapsing the $\sss$-submodule $T(E)$ of the $\sss$-module $E$ according to the following general process. We let $\crel$  be  the category of pairs $X\supset Y$ of pointed sets (here $*\in Y\subset X$) where morphisms are maps of pointed sets $f:X\to X'$ such that $f(Y)\subset Y'$. The collapsing of $Y$ to the base point $*$ defines a functor $\cC:\crel\longrightarrow \Ses$, $\cC(X,Y)= X/Y$.
\begin{defn}\label{scat}
Let $F$ be an $\sss$-submodule of an $\sss$-module $E$. Then the  $\sss$-module $E/F$ is the composition of the pair $(E,F)$ with the functor $\cC$: $E/F(k_+):=E(k_+)/F(k_+)$, $\forall k\in \N$.
\end{defn}
  Next, we notice that each term of the sequence \eqref{exse} is globally invariant  for the action of the endomorphism $T$ and the  induced action  on the right term of the sequence is null \ie its range is reduced to  $*$. Since one divides as in \cite{Al1} by null endomorphisms, this means that the right term of \eqref{exse} is null. Therefore, in this context, it is  natural to define additive invariants  as follows

\begin{defn} \label{defninv} An additive invariant on a class $\cC$ of $\sss$-modules stable under the operation defined in Lemma \ref{scat}, is a map $\chi: \text{End}(\cC)\to G$ from endomorphisms of objects of $\cC$ to an abelian group $G$  that  satisfies the following conditions
\begin{enumerate}
\item $\chi(E,T)=\chi(T(E),T)$
\item $\chi((E_1\vee E_2, T_1\vee T_2)=\chi(E_1,T_1)+\chi(E_2,T_2)$.
\end{enumerate}
\end{defn}
 In this paper we consider the simplest of such class $\cC$, namely the collection   of $\sss$-modules of the form $\sss[F]:=\sss\wedge F$ where $F$ is a  finite pointed set. Any such $\sss$-module is a finite sum (under the join operation $\vee$) of copies of $\sss$.

  We denote by $\W_0(\sss)$ the abelian group {\em range} of the {\em universal} additive invariant  $\tau(E,T)$ for this class $\cC$ of $\sss$-modules. The smash product in $\text{End}(\cC)$ is defined as follows
 $$
 (E_1,T_1)\wedge (E_2,T_2):=(E_1\wedge E_2,T_1\wedge T_2)
 $$
(where the $\wedge$ on the right of the above formula denotes the external smash).  One derives a commutative ring structure on  $\W_0(\sss)$.\newline
 Given an endomorphism $T: F\to F$ of a finite pointed set, we denote by $\tr(T)$ (the trace of $T$) the number of fixed points of $T$ minus one (that accounts for discarding the base point as a trivial fixed point). The subsets $T^n(F)\subseteq F$ obtained by iterating the endomorphism $T$ $n$-times, define a decreasing filtration of $F$ which  stabilizes, after finitely many steps, to $T^\infty(F)\subseteq F$.\newline
  Notice that for  a permutation $T: F\to F$ in $\cC$, the map $\N\ni n\mapsto T^n$ from non-negative integers to iterates of $T$ extends by periodicity to the profinite completion of the integers $\hat \Z:=\varprojlim_n \Z/n\Z$. In this case, we shall use freely the symbol $n\mapsto T^n$ for $n\in \hat\Z$. In the following,  we identify the compact abelian group $\hat\Z$ with  the dual of the discrete abelian  group $\Q/\Z$. More specifically, we first use the ring isomorphism $\hat\Z \simeq \End(\Q/\Z)$ mapping  an integer  $n\in \Z \subset \hat\Z$  to the endomorphism of multiplication by $n$ in  $\Q/\Z$.  Then, we associate to  $n\in  \hat\Z \simeq \End(\Q/\Z)$ the character of the discrete group $\Q/\Z$ obtained by composition with the canonical character $\exp(2\pi i \bullet): \Q/\Z\to S^1$.  One derives the duality pairing
\begin{equation}\label{exse0}
  	\langle s,n\rangle :=\exp(2\pi i ns)\qqq s\in \Q/\Z, \ n\in \hat\Z \simeq \End(\Q/\Z).
  \end{equation}

 Next result determines the ring $\W_0(\sss)$ for the chosen class   $\cC$ of $\sss$-modules of the form $E=\sss[F]:=\sss\wedge F$, with $F$ a finite pointed set.

\begin{thm}\label{w0s} The ring $\W_0(\sss)$ is canonically isomorphic to the invariant part
 of the group ring $\Z[\Q/\Z]$ for the action of the  group  $\Aut(\Q/\Z)= \widehat \Z^\times$. \newline
The universal  additive invariant $\tau: \cC \to \Z[\Q/\Z]$ is uniquely determined by the equality of Fourier transform
\begin{equation}\label{exse2}
	\widehat{\tau(E,T)}(n)=\tr(T^n_{\vert T^\infty(E(1_+))})\qqq n\in \hat\Z.
\end{equation}
\end{thm}
\proof  The evaluation $E\mapsto E(1_+)$ determines a functor  $\smod\longrightarrow \Ses$  and when applied to objects of $\cC$  of the form $E=\sss[F]:=\sss\wedge F$, the evaluation gives  finite pointed sets.  This functor transforms joins $E_1\vee E_2$  to joins $F_1\vee F_2$ and  maps the external smash product of objects of $\cC$ into the smash product of pointed sets (within the category of $\sss$-modules one has $(F_1\wedge \sss)\wedge (F_2\wedge \sss)= (F_1\wedge F_2)\wedge \sss$). In order to show that  $\tau$ defines an additive invariant, whose range is inside $\Z[\Q/\Z]$,  we can  work at the level of the finite pointed sets $F= E(1_+)$. The right hand side of  \eqref{exse2} defines, by construction, an additive map to  functions on $\hat\Z$ that  is also multiplicative. The restriction of $T$ to $T^\infty(F)$ is a bijection preserving the base point, and thus admits a canonical decomposition as a sum of cyclic permutations  of the form $C(k):=((\Z/k\Z)_+,x\mapsto x+1)$, where one adds a base point fixed under $T$. Let $e:\Q/\Z\to \Z[\Q/\Z]$ be the canonical map from the additive group $\Q/\Z$ to its group ring, then we claim that \eqref{exse2} holds if one sets
\begin{equation}\label{exse1}
\tau(C(k)):= \sum_{k\gamma=0}e(\gamma) \in \Z[\Q/\Z].
\end{equation}
Indeed, using \eqref{exse0} the Fourier transform $\widehat{\tau(E,T)}(n)$ is given by
$$
\widehat{\tau(E,T)}(n)=\sum_{k\gamma=0}\langle \gamma,n\rangle =\sum_{ks=0}\exp(2\pi i ns)=\begin{cases} k &\mbox{if } n\in k \hat\Z\\
0  & \mbox{otherwise. }\end{cases}
$$
Thus one gets \eqref{exse2} since the trace of the powers of the cyclic permutation of order $k$ is  $k$ if $ n\in k \hat\Z$ and $0$ otherwise. This argument shows that $\tau$ is an additive invariant of $\cC$ with values in $\Z[\Q/\Z]$. In fact, it also proves that $\tau$ is multiplicative, because the Fourier transform is an injective ring homomorphism from the group ring $\Z[\Q/\Z]$ to the algebra $C(\hat \Z)$ of continuous functions on the compact space $\hat \Z$. It remains to show that $\tau$ is injective and surjective.    \newline
Let $E=\sss\wedge F$ be an object of $\cC$ and $T$ an endomorphism of $E$. Lemma 2.1 of \cite{schemeF1} shows that for $X$ an $\sss$-module, the map associating to $x\in \Hom_\sss(\sss,X)$  its value $x(1_+)(1)\in X(1_+)$  at $1\in 1_+$  defines a bijection of sets $\Hom_\sss(\sss,X)\stackrel{\sim}{\to} X(1_+)$. Since $E$ is a finite sum of copies of $\sss$,  any endomorphism of $E$ is thus of the form $\id\wedge T$ where $T:F\to F$ is a map of  finite pointed sets. Hence the  evaluation functor $(E,T)\mapsto (E(1_+),T_{\vert E(1_+)})$ is faithful. Moreover, the replacement $(E,T)\mapsto (T(E),T)$ corresponds to the replacement $(F,T)\mapsto (T(F),T)$ and by $(i)$ of Definition \ref{defninv} one has the equivalence $(F,T)\simeq (T^\infty(F),T)$. Thus, to show that $\tau$ is injective it is enough to prove that if $(F_j,T_j)$, $j=1,2$, are two permutations of finite pointed sets such that $\tau(F_1,T_1)=\tau(F_2,T_2)$ then they are isomorphic. In fact, it suffices to show that the number $m(k)$ of cycles of length $k$ is the same in both permutations for any integer $k$.  Observe that the number $m(k)$ is given in terms of  $\widehat{\tau(E,T)}(n)$ and the M\"obius function $\mu$ by the equality (that follows by inverting the relation $\sum_{u\vert n} m(u)=\widehat{\tau(E,T)}(n)$)
$$
m(k)= \sum_{d\vert k}\mu(k/d) \widehat{\tau(E,T)}(d).
$$
We have thus shown that the number $m(k)$ is determined by $\tau(E,T)$ and hence, as explained above, that $\tau$ is injective. To show the surjectivity  of $\tau$, it is enough to check that the generators of the abelian group $\Z[\Q/\Z]^G$,  $G=\Aut(\Q/\Z)$, are in the range of $\tau$. It is well known that a system of generators of that group is given by the elements $\rho(n):= \sum_{\gamma \in P(n)}e(\gamma)$, where $P(n)\subset \{\gamma\in \Q/\Z\mid n\gamma=0\}$ is the subset  of primitive roots.   In view of the relation
$$
\rho(n)=\sum_{d\vert k}\mu(k/d) \tau(C(d)),
$$
obtained by inverting the relation
$$
\sum_{u\vert n} \rho(u)=\tau(C(n)),
$$
 one derives that $\tau$ is surjective. \endproof

 Formula \eqref{exse1}  expresses the fact  that $\tau$ associates to an endomorphism $T$ in $\cC$ its  ``trace" $\tau(T)$ viewed as the sum of its {\em abstract} eigenvalues. The Fourier transform  $\widehat{\tau(E,T)}(n)$ as in \eqref{exse2} plays then the role of generalized ghost component: $\widehat{\tau(E,T)}(n)=\gh_{n}((E,T))$, following the classical description of the global Witt ring.

   The Frobenius endomorphisms and the Verschiebung maps  are meaningful operators  in $\W_0(\sss)$. The Frobenius endomorphism $F_n$ is defined by the equality $F_n((E,T)):=(E,T^n)$.  The Verschiebung map $V_n$ replaces the pair $(E,T)$ by the endomorphism of the sum $\vee^n E$ of $n$ copies of $E$ that cyclicly permutes the terms and uses $T:E\to E$ to pass from the last term to the first. In particular, when applied to a cyclic permutation of order $m$ this process yields a cyclic permutation of order $nm$, as in an odometer, and one gets the following relations
\begin{equation}\label{rabi9}
    V_{nm}=V_n\circ V_m=V_m\circ V_n\,, \ \ F_{nm}=F_n\circ F_m=F_m\circ F_n.
\end{equation}
 Next statement displays several  equations connecting  Frobenius and Verschiebung operators on $\W_0(\sss)$.

\begin{prop} \label{rabi10}   The following relations hold with $m,n\in\N$\vspace{.05in}

$(1)$~$F_n\circ V_n(x)=nx$.\vspace{.05in}

$(2)$~$V_n(F_n(x)y)=xV_n(y)$.\vspace{.05in}

$(3)$~If $(m,n)=1$, ~~$V_m\circ F_n=F_n\circ V_m$.\vspace{.05in}

$(4)$~For $n\in \N$, ~~ $V_n(x)V_n(y)=nV_n(x y)$.
\end{prop}

\proof The relations $(1)$ and $(4)$ are easily verified directly. The other two relations can be deduced from Theorem \ref{w0s}: the point being that a direct proof necessarily makes use of condition $(i)$ in Definition \ref{defninv}. In fact using the invariant $\tau$,  we shall  obtain the equalities
\begin{equation}\label{2rho}
\tau(V_n(x))=\tilde\rho_n(\tau(x))\qqq n\in \N, \ x\in \W_0(\sss),
\end{equation}
 where $\tilde\rho_n$, $n\in\N$, is defined by
\begin{equation}\label{newrhos}
\tilde\rho_n: \Z[\Q/\Z] \to \Z[\Q/\Z], \qquad \tilde\rho_n(e(\gamma))=
\sum_{n\gamma'=\gamma}e(\gamma').
\end{equation}
The action of $\tau$ on $F_n$ is described by
\begin{equation}\label{2sigma}
\tau(F_n(x))=\sigma_n(\tau(x))\qqq n\in \N, \ x\in \W_0(\sss),
\end{equation}
where for each $n\in\N$ the group ring endomorphisms $\sigma_n$ is given by
\begin{equation}\label{defnsig}
\sigma_n: \Z[\Q/\Z] \to \Z[\Q/\Z],\qquad \sigma_n(e(\gamma))=e(n\gamma).
\end{equation}
To prove \eqref{2rho} and \eqref{2sigma} we claim that in terms of the Fourier transforms $\gh_{n}((E,T))=\widehat{\tau(E,T)}(n)$ as in \eqref{exse2} one  has
\begin{equation}\label{rhoinf}
\gh_n(V_m((E,T)))=\begin{cases} m \,\gh_{n/m}((E,T)) & \text{if $m\vert n$}\\ 0 &\text{otherwise}\end{cases}
\end{equation}
and
\begin{equation}\label{sigmainf}
\gh_n(\sigma_m((E,T)))=\gh_{nm}((E,T)) .
\end{equation}
Indeed, one verifies \eqref{rhoinf} and \eqref{sigmainf} directly using \eqref{exse2} and also  easily checks that
$
\widehat{\sigma_m(f)}(n)=\widehat{f}(nm)
$.
The Fourier transform of $\tilde\rho_m(e(\gamma))=\sum_{m\gamma'=\gamma}e(\gamma')$ is
$$
\sum_{m\gamma'=\gamma} \exp(2\pi i n \gamma')=\begin{cases} m \,\exp(2\pi i \gamma n/m ) & \text{if $m\vert n$}\\ 0 &\text{otherwise.}\end{cases}
$$
This proves \eqref{2rho} and \eqref{2sigma}. The four equalities of Proposition \ref{rabi10} then follow from Proposition 5.1 of
 \cite{ArithBC} that  is also reported here below for completeness. \endproof

\begin{prop}\label{mapstilderho} The endomorphisms $\sigma_n$ and the maps $\tilde\rho_m$  fulfill the following relations
\begin{equation}\label{rhomult}
 \sigma_{nm}=\sigma_{n}\sigma_{m}\,,\ \
    \tilde\rho_{mn}=\tilde\rho_{m}\tilde\rho_{n}\qqq m,n\in\N
\end{equation}
\begin{equation}\label{rhosigmamult}
    \tilde\rho_{m}(\sigma_m(x)y)=x\tilde\rho_{m}(y)\qqq x,y\in \Z[\Q/\Z]
\end{equation}
\begin{equation}\label{divisible}
    \sigma_c(\tilde\rho_b(x))=(b,c)\,\tilde\rho_{b'}(\sigma_{c'}(x))\,, \ \
    b'=b/(b,c)\,, \ \ c'=c/(b,c)\,,\qquad (b,c):=\text{gcd}(b,c).
\end{equation}
\end{prop}

We note that taking $b=c=n$ in \eqref{divisible} gives
\begin{equation}\label{divisiblebis}
    \sigma_n(\tilde\rho_n(x))=n\,x\qqq x\in \Z[\Q/\Z].
\end{equation}
On the other hand, if one takes $b=n$ and $c=m$ to be relatively prime one has
\begin{equation}\label{divisibleter}
    \sigma_n\circ\tilde\rho_m= \tilde\rho_m\circ \sigma_n.
\end{equation}

\subsection{$\W_0(\sss)$ and the integral BC-system}\label{sectbcsys}

We recall  from \cite{ArithBC}  the  main structure of the integral BC-system.  The integral $BC$-algebra
 is the algebra $\cH_\Z=\Z[\Q/\Z]\rtimes_{\tilde\rho}\N$ generated by the group
ring $\Z[\Q/\Z]$, and by the elements $\tilde\mu_n$ and $\mu_n^*$, with
$n\in\N$,
 which satisfy
  the
relations
\begin{equation}\label{presoverZ2}
\begin{array}{l}
\tilde\mu_{nm}= \tilde\mu_n  \tilde\mu_m \qqq n,m\in\N\\[3mm]
\mu_{nm}^* =\mu_{n}^*\mu_{m}^* \qqq n,m\\[3mm]
\mu_n^* \tilde\mu_n  =n \\[3mm]
\tilde\mu_n\mu_m^* = \mu_m^*\tilde\mu_n \ \ \ \ (n,m) = 1.
\end{array}
\end{equation}
 as well as the relations
\begin{equation}\label{presoverZ1}
\begin{array}{l}
\tilde\mu_n x \mu_n^* = \tilde\rho_n(x)\ \ \ \  \\[3mm]
\mu_n^* x = \sigma_n(x) \mu_n^*  \\[3mm]
x \tilde\mu_n = \tilde\mu_n \sigma_n(x),
\end{array}
\end{equation}
 Notice that the  equations \eqref{rhomult}, \eqref{rhosigmamult}, \eqref{divisible} holding on $\cH_\Z$ are the same as those fulfilled by the Frobenius endomorphisms $F_n$ and the Verschiebung maps $V_n$ on $\W_0(\sss)$: see Proposition~\ref{rabi10}. More precisely, under the correspondences $\sigma_n\to F_n$,
$\tilde\rho_n\to V_n$ the two  equations of \eqref{rabi9} correspond to \eqref{rhomult}, and the first three equations of Proposition \ref{rabi10} correspond respectively to \eqref{divisiblebis}, \eqref{rhosigmamult} and \eqref{divisibleter}. These results evidently point out to the existence of a strong relation between the ($\lambda$)-ring $\W_0(\sss)$ and the  ring $\Z[\Q/\Z]^G$ endowed with the aforementioned operators. This result is in agreement with the idea developed in \cite{ArithBC}  that the integral BC-system is the Witt ring of the algebraic closure of $\F_1$.

\subsection{$\lambda$-rings and the Arithmetic Site} \label{sectlambda}

We refer to \cite{Atiyah1} for the notion of  $\lambda$-ring: here   we are only concerned with  the case of special $\lambda$-rings in the terminology of \opcit  We recall that the Adams operations $\psi_k$ are endomorphisms  defined explicitly in terms of the $\lambda$-operations $\lambda_j$. Thus for instance  we have
$$
\begin{array}{ccc}
 \psi_1 & =& \lambda_1\\
  \psi_2&=&\lambda_1^2-2 \lambda_2  \\
 \psi_3 &=& \lambda_1^3-3 \lambda_2 \lambda_1+3 \lambda_3\\
 \psi_4 &=& \lambda_1^4-4 \lambda_2 \lambda_1^2+4 \lambda_3 \lambda_1 +2 \lambda_2^2-4 \lambda_4\\
\end{array}
$$
To  obtain the general formula one uses (see \opcit \S 5) the  identity of formal power series
\begin{equation}\label{lognewton}
	t\partial_t(\lambda(t))/\lambda(t)=\sum_1^\infty (-1)^{n}\psi_nt^n, \ \
	\lambda(t):=\sum_0^\infty \lambda_nt^n.
\end{equation}
 Next, we observe that  $\lambda$-rings   belong naturally to the topos $\wnt$  that underlies the Arithmetic Site.
\begin{lem}\label{lambda1} Let $R$ be a $\lambda$-ring. Then the Adams operations $\psi_n$ turn $R$ into a sheaf of rings on the topos $\wnt$.
\end{lem}
\proof  The semigroup $\nt$ acts by endomorphisms on $R$ as follows (\hspace{-4pt}\cite{Atiyah1} Proposition 5.2):
$$
(n,x)\mapsto \psi_n(x)\qqq x\in R, \ n \in \nt,
$$
so that $R$, endowed with this action,   becomes  a sheaf of rings on the (presheaf) topos $\wnt$.\endproof
 The definition of the Arithmetic Site  was in fact motivated by the $\lambda$-operations in cyclic theory \cite{Loday} and their role in the conceptual understanding of  Serre's local factors of geometric L-functions at the archimedean places  \cite{CC6}.  The whole construction relied on the cyclic category \cite{CoExt}  and its epicyclic refinement \cite{topos}.  It is a surprising fact that while the classifying spaces of these small categories are non-trivial \cite{CoExt}, the classifying space of the simplicial category is indeed trivial even though the  associated cohomology theory, the Hochschild cohomology, is non-trivial. It is  exactly the theory of Grothendieck topoi   that provides the best geometric framework to understand Hochschild and  cyclic homologies together with  the $\lambda$-operations. This is  achieved using the topos associated to the cyclic category  and its epicyclic refinement.
 The epicyclic topos $\epi$ is  defined by  considering the presheaf topos associated to the {\em opposite} of the epicyclic category $\lbt$. In other words, this is the topos of {\em covariant} functors from $\lbt$ to $\Se$. This choice   is dictated by the following natural construction. A commutative ring $R$ determines a  {\em covariant} functor $R^\#:\fin\longrightarrow\Ab$ from the category of finite sets to  abelian groups. This functor assigns to a finite set $J$ the tensor power $R^{\otimes J}=\bigotimes_{j\in J}R$ and to a map $\phi:J\to J'$ the morphism of abelian groups
  $$
 R^\#(\phi):\bigotimes_{j\in J}R\to \bigotimes_{k\in J'}R, \qquad
 R^\#(\phi)\left(\otimes_{j\in J} x_j \right)=\otimes_{k\in J'}y_k, \ \ y_k=\prod_{\phi(j)=k} x_j.
 $$
   In \cite{CCproj}  we explained in geometric terms  that there is  a natural {\em covariant} functor  $\cP:\lbt\longrightarrow\fin$.  In this construction, the involved geometry is projective geometry in characteristic $1$ and $\lbt$ is interpreted as a category of projective spaces with geometric morphisms. One lets $\zmax$ be the only semifield whose multiplicative group is infinite cyclic.  One shows that the category $\lbt$ is isomorphic to the category whose objects are the  $\zmax$-modules  $\cE_n$ whose underlying additive semigroup is $\zmax$ while the action of $\zmax$ is given by composing the multiplication with the Frobenius endomorphism $\fr_n(x)=x^n$. The morphisms in the category are the projective transformations \cite{CCproj}. Then   the functor $\cP$  associates to a projective space its finite set of points, it is surjective on objects (up to isomorphism) and morphisms.  Then, for any commutative ring $R$,  the composite  $R^\#\circ \cP$  provides   a {\em covariant} functor $R^\natural=R^\#\circ \cP :\lbt\longrightarrow\Ab$ \ie   a sheaf   of abelian groups over the topos $\epi$.  Both the cyclic homology of $R$ and its $\lambda$-operations are completely encoded by the associated sheaf $R^\natural$. \newline
   When  one works with the category $\fin_*$ of finite {\em pointed} sets, instead of $\fin$, the  corresponding structure involves,  in place of $\lbt$  the semi-direct product category $\dop\ltimes \N^\times$. Here we view $\dop$ as the category of finite intervals $[n]^*:=\{0, \ldots, n+1\}$ with $n+2$ elements (\ie $n$ elements besides the smallest and largest) with morphisms preserving the smallest and largest elements.  The action $\sdd_k$ of $k\in\N^\times$ on $\dop$ is understood as the concatenation of $k$ copies of the interval (see \opcit Lemma 2.2), it is the dual of the edgewise subdivision. One has
 a natural covariant functor $\cF:\dop\longrightarrow \gop$ that associates to  an interval $I=\{b,\ldots,t\}$ the pointed set $I_*=I/\sim$ with base point the  class of $b\sim t$ obtained by identifying the smallest and largest  elements of $I$. To any morphism of intervals $f:I\to J$ corresponds the quotient map $f_*$ which preserves the base point. The semi-direct product category $\dop\ltimes \N^\times$ has the same objects as $\dop$ and new morphisms $\pi_n^k:\sdd_k([n]^*)\to [n]^*$.  The functor $\cF$ can be extended to $\dop\ltimes \N^\times$ as follows (see \opcit Proposition 2.8)
\begin{prop} For any $n\geq 0,k\in\N^\times$, let $(\pi_n^k)_*:\cF(\sdd_k([n]^*))\to\cF([n]^*)$ be given by the residue modulo $n+1$. Then the  extension of the functor  $\cF$  on morphisms given by
$$
\phi=\pi_n^k\circ \alpha\mapsto \phi_*:=(\pi_n^k)_*\circ \alpha_*
$$
 determines a functor $\cF:\dop\ltimes \N^\times\longrightarrow \gop$ which is surjective on objects and morphisms.
\end{prop}
The construction of the $\lambda$-operations on Hochschild homology  of a commutative ring \cite{Loday}, only uses the functor $R^\#\circ \cF:\dop\ltimes \N^\times\longrightarrow \Ab$.
 In view of the fact that the presheaf topos $\Gamma^\wedge$ of covariant functors $\gop\longrightarrow \Se$ plays a key role in defining the category $\smod$ of $\sss$-modules, there is the need  to clarify the relations between the various candidates for the geometric counterpart of the ``absolute point",  given by  the following topoi:
 \begin{enumerate}
 \item The Arithmetic Site $\wnt$.
 \item The topos $\Gamma^\wedge$ of covariant functors $\gop\longrightarrow \Se$.
 \item The topos of covariant functors  $\dop\ltimes \N^\times\longrightarrow \Se$.
 \item The epicyclic topos $\epi$.
\end{enumerate}

\section{Probabilistic interpretation of the $\theta$-dimension}\label{secttheta}

In this section we give a probabilistic interpretation of the $\theta$-dimension  of an Arakelov divisors $D$ on  $\overline{\Spec\Z}$. Given an euclidean lattice $\overline{L}:=(L, \Vert.\Vert)$  its  $\theta$-invariant is defined (see \cite{vdG-S,bost}) by
 \begin{equation}\label{thetainv}
 h_\theta^0(\overline{L}):=\log \sum_{v\in L} \exp(-\pi \Vert v\Vert^2).
  \end{equation}
We refer to \cite{bost} for the extension of  the definition and properties of the $\theta$-invariant to the infinite dimensional case  and its interpretation using the thermodynamic formalism.  Here, we  provide  a different interpretation of the same invariant based on a probabilistic way of counting naively the solutions  involving a one-dimensional {\em complex} vector space as customary in Arakelov geometry.
\subsection{The $\sss$-module $E_\xi$}
 Let $E$ be a	 one dimensional complex vector space. Given a vector $\xi \in E$ we define the $\sss$-module
 \begin{equation}\label{countingz}
 E_\xi(k_+):= \{\phi\in E^k\mid
\sum_{x\in k_+}\Vert \phi(x)\Vert \leq \Vert \xi\Vert\}.
\end{equation}
 The above definition is independent of the choice of the norm $\Vert~\Vert$  on  $E$.
 Given  a rank one discrete subgroup $L\subset E$ of the additive group of $E$ one can roughly  measure the size of a vector $\xi \in E$  by comparing it with the unit of length provided by $L$. This yields the quantity
  \begin{equation}\label{countingLz}
 [\xi/L]:=\# \left( HL\cap E_\xi\right)(1_+),
\end{equation}
where $HL$ denotes the Eilenberg-MacLane $\sss$-module associated to the abelian group  $L$ (see \cite{DGM} Example 2.1.2.1). One has $HL(k_+):=L^k$ and the functoriality is that of $L$-valued divisors on  finite pointed sets. Note that this invariant of the pair $(L\subset E, \xi)$ does not depend  on the choice of  the norm on $E$.

\subsection{The Gaussian measure}
 Next we  recall the well-known characterization of the Gaussian measure
\begin{lem}\label{mcomplex3} Let $\cH$ be a	 one dimensional complex Hilbert space.
There exists a unique probability measure $dp$ on $\cH$ fulfilling the following properties
\begin{itemize}
\item $dp$ is rotation invariant.
\item Two real linear forms which are orthogonal are independent random variables.
\item The mean value of the norm is $\frac 12$.
\end{itemize}
\end{lem}
\proof We can assume that $\cH=\C$ with the absolute value as the norm. The measure  $dp(z):=\exp(-\pi \vert z\vert^2)dx dy$, with $z=x+iy$  fulfills  the required conditions. Note that the first two requirements on $dp$  do not use  the norm and fix the density of the measure to be a Gaussian $a\exp(-\pi a \vert z\vert^2)$. Since the  expectation value of the norm on this  probability measure  is
$$
a\int _{-\infty }^{\infty }\int _{-\infty }^{\infty }\sqrt{ x^2+y^2} \exp \left(-\pi a  x^2-\pi  a y^2\right)dxdy=\frac{1}{2 a^{1/2}}
$$
 the  last condition fixes uniquely the normalization constant $a=1$ and hence the measure.\endproof

\subsection{Probabilistic formula for $\exp(h_\theta^0(D))$}

To an Arakelov divisor $D=\sum u_j\{p_j\}+ u\{\infty\}$ on $\overline{\Spec\Z}$   is  associated the euclidean lattice
 \begin{equation}\label{H0bis}
H^0(D):=\left(H^0(\Spec\Z,\cO(D_\fii)),e^{-u} \vvert~\vvert\right)
\end{equation}
where the invertible sheaf $\cO(D_\fii)$ on $\Spec\Z$ is viewed as a subsheaf of the constant sheaf $\Q$ and the norm $\vvert~\vvert$ is the  restriction of the euclidean norm on $\Q$.
\begin{defn}\label{defncountingLz} Let $D$ be an Arakelov divisor on $\overline{\Spec\Z}$. We let $$H^0_\C(D)=(L\subset E,\vvert~\vvert),$$ where $E$ is the one dimensional Hilbert space $E=\C$ with norm $e^{-u} \vvert~\vvert$ and $L\subset E$ is the rank one discrete subgroup $L=H^0(\Spec\Z,\cO(D_\fii)\subset \Q\subset \C$ of formula \eqref{H0bis}.
\end{defn}
The role of the complexification $H^0_\C(D)=(L\subset E,\vvert~\vvert)$ is exhibited by the following  result
\begin{thm}\label{thmcountingLz}
 Let $dp$ be the unique measure on $E$  defined in  Lemma \ref{mcomplex3}. One has
 \begin{equation}\label{averagec}
\exp(h_\theta^0(D))=\int [\xi/L] dp(\xi).
\end{equation}
\end{thm}
\proof The pair $H^0_\C(D)$ associated to $D$ by Definition \ref{defncountingLz}   is unchanged,  up to isomorphism, if one replaces  $D$ with  $D+(q)$ \ie by adding a principal divisor. Thus we can assume that $D=  u\{\infty\}$, with  $u=\deg(D)\in \R$. Then the pair $H^0_\C(D)=(L\subset E,\vvert~\vvert)$ is the same as the pair $(\Z\subset E,\vvert~\vvert)$, where $E$ is the one dimensional Hilbert space $E=\C$ with norm $\vvert z\vvert:=e^{-u} \vert z\vert$. The measure $dp$ on $E$ as in  Lemma \ref{mcomplex3} is then (since $e^{-u}\frac{1}{2 a^{1/2}}=\frac 12$) $$
dp(z):=a\exp(-\pi a\vert z\vert^2)dx dy,\quad  z=x+iy, \ \ \ a=e^{-2u}.
$$
For $z\in E=\C$ one has $[z/\Z]=1+2 \,\text{IntegerPart}(\vert z\vert)$, thus the right hand side of \eqref{averagec}  becomes
$$
\int [\xi/L] dp(\xi)=a\int (1+2 \,\text{IntegerPart}(\vert z\vert))\exp(-\pi a\vert z\vert^2)dx dy, \ \ \ a=e^{-2u}.
$$
 Let $\chi:=1_{[-1,1]}$ be the characteristic function of the unit interval, one has
$$
\text{IntegerPart}(\vert z\vert)=\sum_{\nt}\chi(\frac {n}{\vert z\vert})
$$
and by applying Fubini's Theorem we obtain
$$
\int [\xi/L] dp(\xi)=1+2a\sum_{\nt}\int \chi(\frac {n}{\vert z\vert})\exp(-\pi a\vert z\vert^2)dx dy.
$$
One has
$$
a\int \chi(\frac {n}{\vert z\vert})\exp(-\pi a\vert z\vert^2)dx dy=
a\int_{\vert z\vert\geq n}\exp(-\pi a\vert z\vert^2)dx dy=
2\pi a\int_n^\infty  \exp(-\pi a\rho^2)\rho d\rho=\exp(-\pi a n^2)
$$
so that
$$
\int [\xi/L] dp(\xi)=
1+2\sum_{\nt}\exp(-\pi a n^2).
$$
This gives the required formula using \eqref{thetainv} since, as $ a=e^{-2u}$, $a n^2$ is the square of the norm of the integer $n$ using \eqref{H0bis}.\endproof

\section{$H^1(\overline{\Spec\Z},\cO(D))$ and the Dold-Kan correspondence for $\sss$-modules}\label{sectdks}

In this section we apply the analogue of the Dold-Kan correspondence in the framework of $\sss$-modules to construct a $\Gamma$-space $\bfh(D)$  associated to an Arakelov divisor $D$ on $\overline{\Spec\Z}$.  $\bfh(D)$ is the homological incarnation of the Riemann-Roch problem for  $D$ and  is built  from a natural morphism of $\sss$-modules  canonically  associated to the divisor $D$, as explained in  \S \ref{sectmorphi}.

 \subsection{The structure sheaf of $\overline{\Spec\Z}$}\label{sectspecz}

 We briefly recall the definition  of the structure sheaf of $\sss$-algebras on the Arakelov compactification   $\Spec\Z\subset\overline{\Spec\Z}$  as introduced in \cite{CCprel}. Here, we  also  show that the construction extends to a one parameter deformation by  implementing a real parameter $0<\alpha\leq 1$ in the key proposition  of \opcit that provides the extension of the structure sheaf of $\Spec\Z$ at the  real archimedean place  and introduces the module $\cO(D)$ attached to a divisor $D$.
 \begin{prop}\label{sssalg2} Let $0<\alpha\leq 1$ be a real number and let $R$ be a semiring.\newline
  $(i)$~Let  $\vvert~\vvert$  be  a sub-multiplicative seminorm on $R$. Then $HR$ is naturally endowed with a structure of $\sss$-subalgebra $\vvert HR\vvert_1^\alpha\subset HR$ defined as follows
\begin{equation}\label{subvert}
\vvert HR\vvert_1^\alpha: \Gamma^{{\rm op}}\longrightarrow \Se_*\qquad \vvert HR\vvert_1^\alpha(F):=\{\phi\in HR(F)\mid \sum_{F\setminus \{*\}} \vvert\phi(x)\vvert^\alpha\leq 1\}.
\end{equation}
$(ii)$~Let $E$ be an $R$-semimodule and $\vvert~\vvert$ a seminorm on $E$ such that $\vvert a\xi\vvert\leq \vvert a\vvert \vvert \xi\vvert$, $\forall a\in R$, $\forall \xi \in E$. Then for any $\lambda\in \R_+$ the following formula defines a module $\vvert HE\vvert^\alpha_\lambda$ over $\vvert HR\vvert_1^\alpha$
\begin{equation}\label{subvert1}
\vvert HE\vvert^\alpha_\lambda: \Gamma^{{\rm op}}\longrightarrow \Se_*\qquad \vvert HE\vvert^\alpha_\lambda(F):=\{\phi\in HE(F)\mid \sum_{F\setminus \{*\}} \vvert\phi(x)\vvert^\alpha \leq \lambda\}.
\end{equation}
\end{prop}
 \proof  One needs to show that the maps
$
HE(f)(\phi)(y):=\sum_{x\mid f(x)=y} \phi(x)
$
are compatible with the conditions
$
\sum_{X\setminus \{*\}} \vvert\phi(x)\vvert^\alpha\leq \lambda.
$
For $\alpha\leq 1$  and for any positive real numbers $a,b$  one has the inequality $(a+b)^\alpha\leq a^\alpha+b^\alpha$. Thus, since  the map $a\mapsto a^\alpha$ is increasing, we  obtain, using  the inequality
$$
\vvert \sum_{x\mid f(x)=y} \phi(x)\vvert \leq
\sum_{x\mid f(x)=y} \vvert\phi(x)\vvert
$$
 that
$$
\vvert \sum_{x\mid f(x)=y} \phi(x)\vvert^\alpha \leq \left(\sum_{x\mid f(x)=y} \vvert\phi(x)\vvert \right)^\alpha \leq
\sum_{x\mid f(x)=y} \vvert\phi(x)\vvert^\alpha.
$$
One  obtains, as required
$$
\sum_{Y\setminus \{*\}} \vvert HE(f)(\phi)(y)\vvert^\alpha=\sum_{Y\setminus \{*\}}\vvert \sum_{x\mid f(x)=y} \phi(x)\vvert^\alpha \leq
\sum_{X\setminus \{*\}} \vvert\phi(x)\vvert^\alpha   \leq \lambda.
$$
The sub-multiplicativity of the norm entails the compatibility with the product.
\endproof
\begin{rem}
	For $\alpha >1$ it is no longer true that \eqref{subvert1} defines a sub $\Gamma$-set of $HE$, since taking $\phi(x)=\xi\neq 0$ for all $x$, one has  when $\#\{x\mid f(x)=y\}=n>1$,
$$
\vvert \sum_{x\mid f(x)=y} \phi(x)\vvert^\alpha=n^\alpha \vvert \xi\vvert^\alpha >
\sum_{x\mid f(x)=y} \vvert\phi(x)\vvert^\alpha
=n \vvert \xi\vvert^\alpha.
$$
\end{rem}

 \subsection{The morphism $\phi$}\label{sectmorphi}

 For a function field $\K$  with field of constants $k$, A. Weil gave a proof of the Riemann-Roch formula using  the topological ring of  adeles and Pontrjagin duality. The proof provides, in particular,  an adelic interpretation of $H^1(D)$, for divisors  $D$ of $\K$. For  every place $v\in \Sigma_\K$ of $\K$, the local field $\K_v$ admits a canonical ultrametric norm $\vert ~\vert_v: \K_v^\times\to \R_+^*$:  we let ${\rm Mod}(\K_v)$  the range of $\vert ~\vert_v$. Then, one can view  a divisor  of $\K$  as a map  $D:\Sigma_\K\to \R_+^*$  on the set $\Sigma_\K$ of places of $\K$  such that $D(v)\in {\rm Mod}(\K_v)$, $\forall v \in \Sigma_\K$ and  fulfilling the condition  $D(v)=1$  for all but  finitely many places. A divisor $D$ defines   a compact open subgroup $\cO(D)$ of the additive group $\A_\K$ of adeles of $\K$ by setting
 \begin{equation}\label{weil1}
\cO(D):=\{a=(a_v)\in \A_\K \mid \vert a_v\vert_v \leq D(v),~\forall v \in \Sigma_\K\}.
\end{equation}
 By implementing the diagonal embedding of $\K$ in  $\A_\K$,  one defines the morphism
\begin{equation}\label{weil2}
\psi: \K\oplus \cO(D)\to \A_\K, \ \psi(x,y):= x+y \in \A_\K\qqq x\in \K, y\in \cO(D).
\end{equation}
 The kernel and cokernel of $\psi$ define  respectively the following two natural $k$-vector spaces
 \begin{equation}\label{weil3}
H^0(D):=\cO(D)\cap \K, \ \ H^1(D):=\A_\K/(\K+\cO(D)).
\end{equation}
 These definitions use in a crucial way the ultrametric property of the local norm $\vert ~\vert_v$ and  the  existence of the  field of constants.\newline
 For number fields, the ultrametric property of the local norm fails at the archimedean places  hence the compact subset of adeles as   in \eqref{weil1} is no longer an additive group.  On the other hand, the map $\psi$  as in \eqref{weil2} still  retains a meaning.
 We consider the simplest case of   the  field $\Q$.  To assign a map  $D:\Sigma_\Q\to \R_+^*$  with $D(v)\in {\rm Mod}(\Q_v)$, $\forall v \in \Sigma_\Q$, and such that $D(v)=1$  for  all but finitely many places, is equivalent to define an Arakelov divisor  $D=\sum a_j\{p_j\}+ a\{\infty\}$  on $\spzb$.  Indeed, one sets $D(p_j):=p_j^{a_j}$, $D(v)=1$ for any finite place $v\neq p_j$ $\forall j$, and $D(\infty)=e^a$.  Then \eqref{weil1} defines a compact subset $\cO(D)\subset \A_\Q$. The analogue of \eqref{weil2} is now the map
 \begin{equation}\label{weil4}
\psi: \Q\times \cO(D)\to \A_\Q, \ \psi(x,y):= x+y \in \A_\Q\qqq x\in \Q, y\in \cO(D).
\end{equation}
We let $j: \A_\Q^{\rm f}\to \A_\Q$, $j(x):= (x,0)$, be the embedding of  finite adeles in adeles. Using the ultrametric property of  the local norms at  the finite places one sees that $\cO(D)^{\rm f}:=\cO(D)\cap j(\A_\Q^{\rm f})$ is a compact subgroup of the additive group $\A_\Q$. Let $G=\Q\times\cO(D)^{\rm f}$: one has $\cO(D)^{\rm f}\cap\Q=\{0\}$ since all elements of $\cO(D)^{\rm f}$ have archimedean component $=0$. Thus, the  restriction  of the map $\psi$ as in \eqref{weil4} to $G$   determines an isomorphism of $G$   with  the  subgroup  $\psi(G)=\Q+\cO(D)^{\rm f}$ of $\A_\Q$. Note that $\psi(G)$  is closed  in $\A_\Q$, since $\Q$ is discrete  (hence closed) and $\cO(D)^{\rm f}$ is compact.
  It follows that $G$ plays no role in the definition of  the kernel and cokernel of  $\psi$ and one can reduce the analysis of  this map  by moding out   $G$  as follows
 \begin{lem}\label{lemredpsi} After quotienting both sides of \eqref{weil4} by $G=\Q\times \cO(D)^{\rm f}$, the map $\psi$  becomes
 	 \begin{equation}\label{weil5}
\phi: [-e^a,e^a]\to \R/L, \quad  \phi(x):=xL \in \R/L\qqq x\in [-e^a,e^a]\subset \R,
\end{equation}
where $L\subset \Q\subset \R$ is the lattice
\begin{equation}\label{weil6}
 L= H^0(\Spec\Z,\cO(D_\fii)):=\{q \in \Q\mid \vert q\vert_v \leq D_v \qqq v\neq \infty \}.
\end{equation}
 \end{lem}
\proof We identify (set-theoretically) the adeles $\A_\Q$ with the product $\A^{\rm f}_\Q\times \R$  endowed with projection morphisms   $p_{\rm f}: \A_\Q\to \A^{\rm f}_\Q$ and $p_\infty: \A_\Q\to \R$. The subgroup $p_{\rm f}(\Q)\subset \A^{\rm f}_\Q$ is dense and the subgroup $p_{\rm f}(\cO(D)^{\rm f})\subset \A^{\rm f}_\Q$ is open.  Hence   $p_{\rm f}(\psi(G))=\A^{\rm f}_\Q$.  The kernel of the restriction of $p_{\rm f}\circ \psi$ to $G$ is the group of pairs $(q,a)\in \Q\times \cO(D)^{\rm f}$ such that $p_{\rm f}(q)+p_{\rm f}(a)=0$. Such pairs are determined by the value of $q$ and this  detects the lattice
$$
\{q \in \Q\mid \vert q\vert_v \leq D_v \qqq v\neq \infty \}=L.
$$
  Thus  one  obtains  an isomorphism  of groups   $p:\A_\Q/\psi(G)\stackrel{\sim}{\to} \R/L$  by taking the restriction of the projection $p_\infty$ to the kernel of $p_{\rm f}$.  By construction  one has $\Q\times \cO(D)\simeq G\times [-e^a,e^a]$,  where the quotient by $G$ is identified with $\{0\}\times  [-e^a,e^a]$, while the map $\psi$ of \eqref{weil4} is  replaced by the map $\phi $ of \eqref{weil5}.\endproof

We can  upgrade  the morphism $\phi$ to a morphism of $\sss$-modules.
Let $L\subset \R$ be a lattice, $q:\R\to \R/L$ the quotient morphism of abelian groups,  and $Hq:H\R\to H( \R/L)$ the associated morphism of $\sss$-modules. For $\lambda >0$, and using the notations as in \S\ref{sectspecz}, one has a natural inclusion of $\sss$-modules $\Vert H\R\Vert_\lambda \stackrel{\iota}{\to} H\R$.  Thus  one derives the following morphism of $\sss$-modules which  coincides with the above $\phi$ at level $1$
\begin{equation}\label{cohom}
\phi: \Vert H\R\Vert_\lambda\to H( \R/L), \ \ \phi := Hq\circ \iota.
\end{equation}

\subsection{Dold-Kan correspondence for $\phi:A\to B$}\label{DK}

Let $\phi:A\to B$ be a morphism   of abelian groups. We describe the simplicial abelian group $\cA$ associated,  by the Dold-Kan correspondence, to the following  complex  $\mathfrak C=\{C_n,\phi_n\}$ indexed in  non-negative degrees,  and defined as follows
\begin{equation}\label{short}
\mathfrak C=\{C_1\stackrel{\phi_1}{\to} C_0\};\qquad C_0=B, \ C_1=A, \quad \phi_1(x):=\phi(x)\qqq x\in A.\end{equation}
 The goal  of this section is to obtain a description of $\cA$ that  is easy to transpose when  the above morphism $\phi$ is replaced by the morphism of $\sss$-modules  as in \eqref{cohom}. \newline
We use the following  definition of the simplicial abelian group $\cA$ associated, by the  Dold-Kan correspondence, to  a complex $C$ of abelian groups
$$
\cA_n:=\Hom_{{\rm Ch}_+}\left(N(H\Z\circ \Delta^n),C\right) \qquad\forall n\ge 0,
$$
where $H\Z\circ \Delta^n$ is the free simplicial abelian group associated to the $n$-simplex (\ie the simplicial set $\Hom(\bullet, [n])$),  and  $N$  denotes the {\em normalized} Moore complex of a simplicial abelian group.  There is  a concrete description of $\cA_n$ given as follows (see \cite{GJ} \S III.2, Prop. 2.2)

\begin{equation}\label{descript}
\cA_n=\bigoplus_{\cF(n,k)} C_k, \quad \cF(n,k):=\{\sigma\in \Hom_\Delta([n],[k])\mid \sigma([n])=[k]\},
\end{equation}
here the direct sum repeats the term $C_k$ of the chain complex as many times as  the number of elements of  the set $\cF(n,k)$ of surjective morphisms $\sigma\in \Hom_\Delta([n],[k])$. In the simple case of the short complex  $\mathfrak C$ in \eqref{short}, the allowed values of $k$ are $k=0,1$.  In this case, the set $\cF(n,0)$ is reduced to a single element:  $\Hom_\Delta([n],[0])=:\Delta_n^0$.  Next, we determine $\cF(n,1)$.
 An element  $\xi\in\Hom_\Delta([n],[1])=\Delta^1_n$ is characterized by  $\xi^{-1}(\{1\})$  that is an hereditary subset of $[n]$. Thus the  vertices in $\Delta^1_n$  are   labelled as
\begin{equation}\label{identxi}
\xi_j, \ 0\leq j\leq n+1, \ \xi_j(k)=1\iff k\geq  j.
\end{equation}
 For each integer $n\geq 0$  the finite set $F(n):=\cF(n,1)$ of surjective elements  $\xi\in\Delta^1_n$  excludes $\xi_{0}$ and  $\xi_{n+1}$, thus $F(n)=\{\xi_j\mid j\in \{1,\ldots, n\}\}$ has $n$ elements.  In the simple case of the complex $\mathfrak C$ one can rewrite \eqref{descript}   as
 \begin{equation}\label{descriptshort}
\cA_n=HB(1_+)\oplus HA(F(n)_+).
\end{equation}
  For the  description of  the simplicial structure,   namely the definition  for each $\theta \in
 \Hom_\Delta([m],[n])$ of a map of sets $\cA(\theta):\cA_n\to \cA_m$  we follow again  \cite{GJ} \S III.2, p. 160-161.  Let, as in \S \ref{sectw0}, $\crel$  be  the category of pairs  of pointed sets $(X,Y)$, with $X\supset Y$  where morphisms are maps  of pairs of pointed sets, and $\gamma:\crel\longrightarrow \Ses$ the functor $(X,Y)\to X/Y$  of collapsing  $Y$ to the base point. One has  the following
 \begin{prop}\label{relativeH} Let $\phi:A\to B$ be a morphism   of abelian groups. The following setting defines a covariant functor
 \begin{equation}\label{descript1}
  H_\phi:\crel\longrightarrow \Ab,  \qquad	H_\phi(X,Y):=HB(Y)\times HA(X/Y).
 \end{equation}
 For  a morphism 	$f:(X,Y)\to (X',Y')$   of pairs of pointed sets, lets, with $\psi= (\psi_Y,\psi_{X/Y})\in H_\phi(X,Y)$, $y'\neq *$,
 \begin{equation}\label{descript2}
 	H_\phi(f)(\psi):=\left(\psi_{Y'}, HA(f)(\psi_{X/Y})\right), \  \psi_{Y'}(y'):=HB(f)(\psi_{Y})(y')+\sum_{x\in X\setminus Y, f(x)=y'}\phi(\psi_{X/Y}(x)).
 \end{equation}
 \end{prop}
\proof  The sum in  \eqref{descript2} is well defined and determines  an element of $B$ for each $y'\in Y'$, $y'\neq *$. By construction $H_\phi(f)$ is a morphism of abelian groups.  Let $g:(X',Y')\to (X'',Y'')$ be a morphism in $\crel$,  and $h=g\circ f$. Let $y''\in Y''$, $y''\neq *$. We show that $H_\phi(h)=H_\phi(g)\circ H_\phi(f)$. We view the elements $\psi= (\psi_Y,\psi_{X/Y})\in H_\phi(X,Y)$ as functions, with $\psi(x)\in A\cup B$,  for $x\in X$ such that
\begin{equation}\label{descript2.5}
\psi(x)\in A \qqq x\in X\setminus Y, \quad \text{and}\ \ \psi(y)\in B\qqq y\in Y, y\neq *, \ \psi(*)=0.
 \end{equation}
We  define $p:A\cup B\to B$  to be equal to $\phi$ on $A$ and  the identity on $B$. Then one has
\begin{equation}\label{descript3}
 	p(\sum a_j)=\sum p(a_j) \qqq a_j\in A,\ \ p(\sum p(b_j))=\sum p(b_j)\qqq b_j\in A\cup B.
 \end{equation}
  It follows from \eqref{descript2}   that
 \begin{equation}\label{descript4}
H_\phi(f)(\psi)(y)= \begin{cases} \sum_{x\in X, f(x)=y} \psi(x) &\mbox{if } y\notin Y'\\
\sum_{x\in X, f(x)=y} p(\psi(x))  & \mbox{if }y \in Y', y\neq *\end{cases}
 \end{equation}
 By applying the functoriality of $HA$ and of the collapsing functor $(X,Y)\longrightarrow X/Y$, one checks that the equality $H_\phi(h)(\psi)=H_\phi(g)( H_\phi(f)(\psi))$ holds for the evaluation on
$x\in X\setminus Y$. Next, we  take $y\in Y, y\neq *$ and evaluate both sides
$$
H_\phi(h)(\psi)(y)=\sum_{x\in X, h(x)=y} p(\psi(x)) =\sum_{z\in X, g(z)=y}\ \ \sum_{x\in X, f(x)=z} p(\psi(x))
$$
$$
H_\phi(g)( H_\phi(f)(\psi))(y)=\sum_{z\in X, g(z)=y}p(H_\phi(f)(\psi)(z)).
$$
Using \eqref{descript3} and \eqref{descript4} one obtains
$$
p(H_\phi(f)(\psi)(z))=\sum_{x\in X, f(x)=z} p(\psi(x)).
$$
This  argument shows that $H_\phi(h)=H_\phi(g)\circ H_\phi(f)$,  thus  proving the required functoriality of $H_\phi$. \endproof

As already mentioned in  \S\ref{sectlambda}, the  category $\dop$ is naturally isomorphic (hence identified) to (a skeleton of) the category of finite intervals, where  by an interval  we mean  a totally ordered set with a smallest element distinct from a largest element.  A morphism of intervals is a non-decreasing map  that preserves the smallest and largest elements. The canonical contravariant functor $\Delta\longrightarrow \dop$  that identifies the opposite category of  $\Delta$ with  $\dop$ (described by intervals as above)  maps the  finite ordinal  $[n]=\{0,1,\ldots,n\}$  to the interval $[n]^*:=\{0, \ldots, n+1\}$.  The dual of  $\theta \in
 \Hom_\Delta([n],[m])$  is
\begin{equation}\label{ident0}
\theta^* \in  \Hom_\dop([m]^*,[n]^*)), \ \ j\leq \theta(i)\iff \theta^*(j)\leq i.
 \end{equation}
 The bijection  of sets  $ \Hom_\Delta([n],[m])\stackrel{\theta\mapsto \theta^*}{\longrightarrow} \Hom_\dop([m]^*,[n]^*)$ implements  the identification
 \begin{equation}\label{ident1}
 \Delta^1_n= \Hom_\Delta([n],[1])\stackrel{\sim}{\to}\Hom_\dop([1]^*,[n]^*)=[n]^*
 \end{equation}
 and one  easily checks that under this identification the map $\xi_j$ of \eqref{identxi} corresponds to $j\in [n]^*=\{0, \ldots , n+1\}$.

We  let  $\partial:\dop\longrightarrow \crel $ be the functor that  replaces an interval $I$ by the pair $\partial I:=(X,Y)$, where $X$ is the set $I$ pointed by its smallest element, and $Y\subset X$ is the subset formed by the smallest and largest elements of $I$. \newline
 Next proposition describes the simplicial abelian group that corresponds to the complex $\mathfrak C$, by the Dold-Kan correspondence
 \begin{prop}\label{DKsimple} Let $\phi:A\to B$ be a morphism   of abelian groups. Then, the  simplicial abelian group associated by the Dold-Kan correspondence to the short complex $\mathfrak C$ in \eqref{short} is canonically isomorphic to the  composite functor $H_\phi\circ \partial: \dop\longrightarrow \Ab$  with
 	$H_\phi:\crel\longrightarrow \Ab$ as in Proposition \ref{relativeH}.
 	 \end{prop}
\proof  For $n\in \N$, the pair $\partial [n]^*$ is identified, using \eqref{ident1} and \eqref{identxi}, to the pair $(X,Y)$  with
$$
X:=\Delta^1_n= \Hom_\Delta([n],[1]), \ \ Y:= \{\xi_0, \xi_{n+1} \},
$$
with  base point  $\xi_0\in Y$,  so that one  obtains $X\setminus Y=F(n)$. It follows from
 \eqref{descriptshort} and \eqref{descript1} that one has a canonical isomorphism of abelian groups
 \begin{equation}\label{identiff}
 \cA_n=HB(Y)\oplus HA(X/Y)=HB(Y)\times HA(X/Y)=H_\phi(X,Y)= H_\phi(\partial [n]^*).
 \end{equation}
 Let $\theta \in
 \Hom_\Delta([m],[n])$, and $\cA(\theta):\cA_n\to \cA_m$ as described in
 \cite{GJ} \S III.2, p. 160-161.  We show that $\cA(\theta)=H_\phi\circ \partial( \theta^*)$, where $\theta^*$ is  defined in \eqref{ident0}. By construction, $\cA(\theta)$ is a morphism of abelian groups  that is determined by its restriction to each copy $C_{j,\sigma}$ of $C_j$ indexed by a surjection $\sigma:[n]\to [j]$ as follows. One writes $ \sigma\circ \theta$ uniquely as the  composite $d\circ s$
\begin{equation}\label{uniqdecomp}
 \sigma\circ \theta= d\circ s, \ \ s\in\Hom_\Delta([m],[k]), \ d\in\Hom_\Delta([k],[j]),
\end{equation}
 where $s$ is a surjection and $d$ an injection. Then, for $\alpha\in C_{j,\sigma}$,  $\cA(\theta)(\alpha)=\cA(d)(\iota_{s,\sigma}\alpha)$,  where $\iota_{s,\sigma}$ is the identity map of $C_j$  that identifies $C_{j,\sigma}$  with $C_{j,s}$. The map $\cA(d)$ is the identity when $d=\id$ (\ie when $k=j$) and is $0$ otherwise, unless $k=j-1$
 and $d\in\Hom_\Delta([j-1],[j])$ is  the injection which misses $j$. In  the case of the complex $\mathfrak C$, the only  indices are $j=0,1$. For $j=0$, there is only one element $\sigma:[n]\to [j]$,  $\sigma\circ \theta$ is surjective and is the only element of $\Hom_\Delta([m],[0])$. This shows that the restriction of $\cA(\theta)$ to $HB(Y)$ is the identity. For $j=1$, one has
  $\sigma\in F(n)$ \ie  $\sigma:[n]\to [1]$ is surjective. If the same result holds for $\sigma\circ \theta :[m]\to [1]$ then, by reasoning as above, one concludes that $\cA(\theta)$ is the natural identification of $C_{1,\sigma}$  with $C_{1,\sigma\circ \theta}$.   When $\sigma\circ \theta :[m]\to [1]$ fails to be surjective there are two cases to consider. One has either $\sigma\circ \theta=\xi_{0}$ or $\sigma\circ \theta=\xi_{m+1}$. The only case in which the obtained map $\cA(\theta)$ is non-zero on $C_{1,\sigma}$ is when the morphism $d$  in the decomposition $\sigma\circ \theta= d\circ s$ of \eqref{uniqdecomp} is $\delta_1 \in\Hom_\Delta([0],[1])$,  with $\delta_1$  missing $1$. This happens when  $\sigma\circ \theta=\xi_{m+1}$ since it corresponds to a mono $\circ$ epi decomposition of the form
$$
\sigma\circ \theta=\delta_1\circ \tau, \ \ \delta_1 \in\Hom_\Delta([0],[1]), \ \ \tau\in\Hom_\Delta([m],[0].
$$
In that case,  $\cA(\theta)$ restricted to  $C_{1,\sigma}$ is the map to  $C_{0}$ given by the boundary $\phi:C_1\to C_0$ of the complex $\mathfrak C$.  By adding together these various contributions of $\cA(\theta)$  one derives the following formula (using the identification \eqref{identiff} and the map $p:A\cup B\to B$  equal to $\phi$ on $A$ and the identity on $B$)
\begin{equation}\label{descript6}
\cA(\theta)(\psi)(y)= \begin{cases} \sum_{x\in X  \mid \, x\circ \theta=y} \psi(x) &\mbox{if } y\notin Y\\
\sum_{x\in X \mid \, x\circ \theta=y} p(\psi(x))  & \mbox{if }y \in Y, y\neq *\end{cases}
 \end{equation}
 Since $\theta^*(x)=x\circ \theta$, one thus obtains $\cA(\theta)=H_\phi\circ \partial( \theta^*)$  as required. \endproof

 The covariant functor $H_\phi: \crel\longrightarrow \Ab$, when composed with the Eilenberg-MacLane functor $H$ defines a functor $HH_\phi :\crel\longrightarrow \smod $   whose explicit description is provided by the following

 \begin{lem}\label{HH} The functor $HH_\phi :\crel\longrightarrow \smod $ is naturally isomorphic to the functor $((U\circ H_\phi)\times \id)\circ \gamma$, where $U:\Ab\longrightarrow \Ses$ is the forgetful functor and $\gamma:\crel \longrightarrow\Gamma\crel$  associates to an object $(X,Y)$ of $\crel$ the covariant functor
 $$
\gamma(X,Y):\gop\longrightarrow \crel,\quad k_+\mapsto (X\wedge k_+,Y\wedge k_+).
 $$
  \end{lem}
\proof  Both $HH_\phi=H\circ H_\phi$ and $((U\circ H_\phi)\times \id)\circ \gamma$ are functors from  $\crel$  to the category  of $\sss$-modules. More precisely,  the functor $((U\circ H_\phi)\times \id)\circ \gamma$, obtained by applying $U\circ H_\phi$ ``pointwise"   takes values in the category  $\Gamma\Ses$  of pointed $\Gamma$-sets  that is the same as the category $\smod$. The functor $H_\phi $ fulfills, for any pointed map $\rho:k_+\to \ell_+$, the following equalities
\begin{equation}\label{HH1new}
H_\phi(X\wedge k_+,Y\wedge k_+)=H_\phi(X,Y)^k, \ \ H_\phi(\id\wedge \rho)=H(H_\phi(X,Y),\rho).
\end{equation}
 Indeed, elements of $H_\phi(X\wedge k_+,Y\wedge k_+)$ are by construction functions $\psi(x,j)\in A \cup B$, for $x\in X\setminus\{*\}$, $j\in \{1,\ldots ,k\}$ with values in $A$ when $x\notin Y$, and in $B$ for $x\in Y$. Equivalently, these elements  are described  by $k$-tuples of functions $\psi_j\in H_\phi(X,Y)$, where $\psi_j(x):= \psi(x,j)$ for $x\in X\setminus\{*\}$ and $j\in \{1,\ldots ,k\}$. The map $\tau=\id\wedge \rho: (X\wedge k_+,Y\wedge k_+)\to (X\wedge \ell_+,Y\wedge \ell_+)$ is a map of pairs  that  preserves the subset  and also its complement. For such maps $\tau$ the formula \eqref{descript2} simplifies to
 $$H_\phi(\tau)=\left(HB(\tau), HA(\tau)\right).$$
 This argument  proves the second equality in \eqref{HH1new}. This equality gives the isomorphism of the $\sss$-modules  $U\circ H_\phi\circ \gamma(X,Y)\simeq HH_\phi(X,Y)$ for a fixed $(X,Y)$. To  complete the proof let us show that this isomorphism is a natural transformation of functors. Let $f:(X,Y)\to (X',Y')$ be a morphism  in $\crel$. One has $\gamma(f)=f\wedge \id$ as a morphism from $(X\wedge k_+,Y\wedge k_+)$ to $(X'\wedge k_+,Y'\wedge k_+)$.  Moreover
 \begin{equation}\label{HH2}
 H_\phi(f\wedge \id)=H_\phi(f)^k:H_\phi(X,Y)^k\to H_\phi(X',Y')^k
\end{equation}
 follows directly from \eqref{descript4} using $(f\wedge \id)(x,j)=(f(x),j)$ $\forall x\in X, j\in \{1, \ldots, k\}$. Finally, \eqref{HH2} shows that the identification \eqref{HH1new}  is a natural transformation of functors. \endproof

 By combining Proposition \ref{DKsimple} with Lemma \ref{HH}  we finally obtain
 \begin{thm}\label{DKsimple1} Let $\phi:A\to B$ be a morphism   of abelian groups. The  $\Gamma$-space associated by the Dold-Kan correspondence to the short complex  $\mathfrak C$ in \eqref{short} is canonically isomorphic to the functor $$((U\circ H_\phi)\times \id)\circ \gamma\circ \partial : \dop \longrightarrow \smod.$$
 \end{thm}
\proof By construction this  $\Gamma$-space is obtained by applying the Eilenberg-MacLane functor $H$ to the simplicial abelian group associated by the Dold-Kan correspondence to the short complex  $\mathfrak C$. Thus by  Proposition \ref{DKsimple} it is the functor $HH_\phi\circ \partial: \dop \longrightarrow \smod$. Lemma \ref{HH} then allows one to replace $HH_\phi$ by $((U\circ H_\phi)\times \id)\circ \gamma$  and this  determines the required equality.\endproof

  As indicated by  Theorem~\ref{DKsimple1},  the construction of the $\Gamma$-space associated  to the short complex $\mathfrak C$ in \eqref{short}  by the Dold-Kan correspondence involves the morphism $\phi:A\to B$ of abelian groups {\em only} through  the composite functor $U\circ H_\phi:\crel \longrightarrow \Ses$. In the next section we show  that the latter functor  is still meaningful  when  $\phi:A\to B$ is replaced by the morphism \eqref{cohom} of $\sss$-modules   as in \S\ref{sectmorphi}.

\subsection{The $\Gamma$-space $\bfh(D)$}\label{sectgammadiv}
 In this section we associate to an Arakelov divisor  $D$ on $\spzb$ a $\Gamma$-space $\bfh(D)$ that embodies the Dold-Kan counterpart of the short complex \eqref{cohom}
associated to the Riemann-Roch problem for the divisor. This construction is based on Theorem \ref{DKsimple1} and the following

 \begin{lem}\label{HH1} Let $L\subset \R$ be a discrete subgroup and $\lambda\in \R_+^*$. Let $\phi:\R\to \R/L$ be the quotient map. The following  equality  defines a subfunctor   of $U\circ H_\phi:\crel \longrightarrow \Ses$
	\begin{equation}\label{simpconstruct1}
K_{L,\lambda}:\crel \longrightarrow \Ses ,\quad K_{L,\lambda}(X,Y):=\{\psi \in H_\phi(X,Y)\mid \sum_{x\in X\setminus Y}\vert \psi (x)\vert \leq \lambda\}.
\end{equation}
 \end{lem}
\proof  The setting \eqref{simpconstruct1} is well defined  in view of \eqref{descript2.5},   showing that
$\psi(x)\in \R, \forall x\in X \setminus Y$. Let $f:(X,Y)\to (X',Y')$ be a morphism  in $\crel$, then by \eqref{descript4} one has
$$
H_\phi(f)(\psi)(y)=  \sum_{x\in X, f(x)=y} \psi(x) \qqq  y\notin Y'.
 $$
The key fact here is that since $f:(X,Y)\to (X',Y')$ is a morphism of pairs, the following implication holds
\begin{equation}\label{simpconstruct2}
y\in X'\setminus Y' \ \& \ f(x)=y \  \Rightarrow x\in X\setminus Y,
\end{equation}
 from which one derives  the inequality
$$
\sum_{y\in X'\setminus Y'}\vert H_\phi(f)(\psi)(y)\vert\leq \sum_{x\in X\setminus Y}\vert \psi (x)\vert.
$$
 This proves  that $K_{L,\lambda}:\crel \longrightarrow \Ses $ is a subfunctor of $U\circ H_\phi$. \endproof

 With the notation of \S\ref{DK} and in particular as in Lemma~\ref{HH}, we have

\begin{prop}\label{propgam} Let $D=\sum a_j\{p_j\}+ a\{\infty\}$ be a divisor on $\spzb$  and let $\lambda=e^a$.  Consider the abelian group
$$ L=H^0(\Spec \Z,\cO(D_\fii))\subset \Q\subset \R.$$
The functor
\begin{equation}\label{propgam0}
\bfh(D):=(K_{L,\lambda}\times \id)\circ \gamma\circ \partial : \dop \longrightarrow \smod
\end{equation}
 defines a $\Gamma$-space  that depends only on the  linear equivalence class of $D$.
 \end{prop}
\proof  The sheaf $\cO(D_\fii)$ on $\Spec\Z$ is viewed as a subsheaf of the constant sheaf $\Q$: this  gives meaning to the inclusion $L=H^0(\Spec \Z,\cO(D_\fii))\subset \Q$.  Let $q\in \Q_+^*$ and $D'=D+(q)$, then  $L'=q^{-1}L$  and $\lambda'=q^{-1}\lambda$. The multiplication by $q^{-1}$  induces an isomorphism of  functors $K_{L,\lambda}\to K_{L',\lambda'}$.   This shows that $\bfh(D)$ only depends on the equivalence class of $D$ modulo principal divisors.\endproof

\subsection{$\bfh(D)$ at level 1}\label{sectlevel1}

When evaluated on the object $1_+$ of $\gop$, the $\Gamma$-space $\bfh(D)$  defines a simplicial set $\bfh(D)(1_+)$. By construction, this  is a sub-simplicial set of the Kan simplicial set associated to the morphism of abelian groups (viewed as a short complex) $\phi: \R\to \R/L$,  however `per-se' it is not a Kan simplicial set. Thus one needs to   exert care  when  considering the homotopy  of $\bfh(D)(1_+)$:  we refer to \cite{CCgromov} \S 2.1,  for the general notions  and notation  used here as follows, when referring to the homotopies. We simply  recall that the new homotopy theory  $\pik_n$ introduced in \opcit is reduced to $\pik_0$ by implementing  the endofunctor $\Omega$ of ``decalage" in  the category $\sses$ of pointed simplicial sets,  as follows
\begin{equation}\label{pindef}
\pik_n(X,\star):=\pik_0(\Omega^n(X)).
\end{equation}
 One obtains, in this way, a concrete description of the elements of $\pik_n(X,\star)$  and of the  following homotopy relation:
\begin{enumerate}
\item	The $0$-skeleton $(\Omega^n(X))_0$ is the set of simplices $x\in X_n$  such that  $\partial_j(x)=*$  $\forall j$.
\item $(\Omega^n(X))_0$ coincides with ${\Hom}_\sses(S^n,X)\subset X_n$, where   $S^n$ is obtained by collapsing the boundary $\partial \Delta[n]$ of the standard simplex to a single base point.
\item The relation $\cR$ on $(\Omega^n(X))_0={\Hom}_\sses(S^n,X)\subset X_n$ given by
 $$
x \cR y \iff \exists z\in (\Omega^n(X))_1~\text{s.t.}~ \partial_0 z=x~\text{and}~ \partial_1 z=y
 $$
coincides with the relation of homotopy between $n$-simplices  as in  the following formula
\begin{equation}\label{Maydefn}
  \cR=\{(x,y)\in X_n\times X_n\mid \partial_j x=\partial_jy\, \forall j \, \&\,  \exists z \mid \partial_jz=s_{n-1}\partial_j x\, \forall j<n, \ \partial_nz=x,\, \partial_{n+1}z=y\}.
  \end{equation}
\end{enumerate}

When  working with  this (new)  homotopies $\pi_n$ usually it  happens that the  relation $\cR$  as in  \eqref{Maydefn} fails to be an equivalence relation, since the pointed simplicial sets fail in general to be fibrant. Even though  (the new)  $\pi_0$ cannot be directly described as a quotient, we  can still  specify its cardinality,  using the following general definition

\begin{defn}\label{cardpi0} Let  $\cR$ be a reflexive relation on a set $E$.  Define $\#(E/\cR)$  to be the maximal cardinality of subsets $F\subset E$ such that $x,y\in F \ \& \ x\neq y\Rightarrow (x,y)\notin \cR$.
\end{defn}

 Obviously, this definition determines the cardinality of the quotient $E/\cR$ when  $\cR$ is an equivalence relation. As an example, consider the relation
$d(x,y)\leq \epsilon$ in a metric space $(E,d)$. Then the number $\#(E/\cR)$ is called the {\em packing number} and is denoted by $\cM(\epsilon)$.   It plays an important role in the notion of $\epsilon$-capacity \cite{KT}.

The following proposition supplies in particular a first approximation on the definition and the cardinality of $H^1(\overline{\Spec\Z},\cO(D))$.

\begin{prop} \label{proppi0}
The simplicial set $\bfh(D)(1_+)$  fulfills the following properties:
\begin{enumerate}
	\item  $\pi_0(\bfh(D)(1_+))=*$ for  $\deg D\geq -\log 2$.
 \item For  $\deg D<0$, $\#(\pi_0(\bfh(D)(1_+)))$ is the largest integer $n< e^{-\deg D}$.
  \item $\pi_1(\bfh(D)(1_+))$	is in canonical bijection with  the global sections $H^0(\spzb ,\cO(D))$.
\item $\pi_n(\bfh(D)(1_+))$ is trivial for $n>1$.
\end{enumerate}
\end{prop}
\proof
(i)~We can assume, by working up-to linear equivalence, that $D_\fii=0$ so that $D=a\{\infty\}$  with $a=\deg D$. Then, $L=H^0(\Spec \Z,\cO(D_\fii))=\Z\subset \R$.  By construction, the pointed simplicial set $\bfh(D)(1_+)$ is obtained by composing the functor $K_{L,\lambda}:\crel \longrightarrow \Ses $ ($\lambda=e^a$) with $\partial:\dop \longrightarrow\crel$.  In degree $0$ and with the notation of \eqref{identxi}, it is described for the  $\xi_j\in \Hom_\Delta([n],[1])$, $0\leq j\leq n+1$, by
$$
K_{L,\lambda}(\partial[0]^*)=K_{L,\lambda}((*,\xi_1),(*,\xi_1))=\R/L.
$$
  In degree $1$  it is given by
\begin{equation}\label{pideg1}
K_{L,\lambda}(\partial[1]^*)=K_{L,\lambda}((*,\xi_1,\xi_2),(*,\xi_2))=\{(\psi_1,\psi_2)\in \R\times \R/L\mid \vert \psi_1\vert \leq \lambda\}.
\end{equation}
One has $\xi_{0}\circ  \delta_j=\xi_{0}$ for all $j=0,1$ and
$$
\xi_1\circ  \delta_0=\xi_0, \ \ \xi_2\circ  \delta_0=\xi_1, \ \ \xi_1\circ  \delta_1=\xi_1, \ \ \xi_2\circ  \delta_1=\xi_1.
$$
Thus, with $\partial_j: K_{L,\lambda}(\partial[1]^*)\to K_{L,\lambda}(\partial[0]^*)$ the corresponding boundaries, one gets
\begin{equation}\label{pi0}
\partial_0(\psi)=\psi_2, \ \ \partial_1(\psi)=\psi_1+\psi_2 \qqq \psi=(\psi_1,\psi_2)\in K_{L,\lambda}(\partial[1]^*).
\end{equation}
This argument shows  that the relation defined on elements of $K_{L,\lambda}(\partial[0]^*)$ by
\begin{equation}\label{pizeroR}
\cR(\alpha,\beta)\iff \exists 	\psi \in K_{L,\lambda}(\partial[1]^*)\mid \partial_0(\psi)=\alpha \  \& \ \partial_1(\psi)=\beta
\end{equation}
is,  in general, not an equivalence relation.
  The relation $\cR$ of \eqref{pizeroR} is  described, using \eqref{pi0} as follows:
$$
(\psi,\psi')\in \cR\subset \R/L\times \R/L\iff  \exists \alpha\in \R, \ \vert \alpha\vert \leq \lambda, \ \psi'=\psi + \alpha.
$$
When $\lambda\geq 1/2$, \ie $\deg D\geq -\log 2$ this relation becomes  the equivalence relation with quotient a single element.\newline
(ii)~With the notations of (i), $\#(\pi_0(\bfh(D)(1_+)))$ is the packing number of the metric space $\R/\Z$ (endowed with the quotient metric of the usual metric on $\R$) for the constant $\lambda$. For any subset fulfilling the separation condition of Definition \ref{cardpi0}, the balls of radius $\lambda/2$ are pairwise disjoint and hence $n \lambda <1$. Conversely, if $n \lambda <1$ the $n$ points
$x_j=\frac{j}{n},~ 0\leq j<n
$
are such that, in the metric space $\R/\Z$, $d(x_i,x_j)>\lambda$ for $i\neq j$.  This  shows (ii).\newline
(iii)~With the above notations, $\bfh(D)(1_+)$ is described in degree $1$ by \eqref{pideg1}  and by \eqref{pi0} one has:
$
\partial_0(\psi)=\psi_2, \ \ \partial_1(\psi)=\psi_1+\psi_2
$.
Thus the spherical conditions $\partial_j(\psi)=*$ mean that $\psi_2=0$ and that the class of $\psi_1$ in $\R/\Z$ is zero \ie  $\psi_1\in L=\Z$.  Hence the spherical elements are in canonical bijection with  the global sections $H^0(\spzb ,\cO(D))$. It remains to describe the relation \eqref{Maydefn} between the spherical elements. With $X_n=K_{L,\lambda}(\partial[n]^*)$,  this is
\begin{equation}\label{Maydefn2}
  \cR=\{(\phi,\phi')\in X_1\times X_1\mid \partial_j \phi=\partial_j\phi'=*\, \forall j \, \&\,  \exists \psi\in X_2 \mid \partial_0 \psi =*,\,   \partial_1 \psi=\phi,\, \partial_{2}\psi=\phi'\}.
  \end{equation}
   $K_{L,\lambda}$ is described in degree $2$ by
\begin{equation}\label{pidegtwo}
K_{L,\lambda}(\partial[2]^*)=K_{L,\lambda}((*,\xi_1,\xi_2,\xi_3),(*,\xi_3))=\{(\psi_1,\psi_2,\psi_3)\in \R\times \R\times \R/L\mid \vert \psi_1\vert +\vert \psi_2\vert\leq \lambda\}.
\end{equation}
 The three boundary maps take values in $\R\times \R/L$ and are given by
$$
\partial_0(\psi)=(\psi_2,\psi_3), \ \partial_1(\psi)=(\psi_1+\psi_2,\psi_3), \ \partial_2(\psi)=(\psi_1,\psi_2+\psi_3).
$$
The condition $\partial_0(\psi)=*$ of \eqref{Maydefn2} implies that $\psi_2 =0$. Hence  the relation coming from the pairs $(\partial_1(\psi),\partial_2(\psi))$ is  the identity. This proves (iii).\newline
(iv)~For $n>1$ one has
\begin{equation}\label{pidegn}
K_{L,\lambda}(\partial[n]^*)=\{(\psi_j)_{j\in \{1,\ldots,n\}}\in \R^{n-1}\times \R/L\mid \sum_{j\in \{1,\ldots,n-1\}} \vert \psi_j\vert \leq \lambda\}.
\end{equation}
 The boundaries $\partial_j:K_{L,\lambda}(\partial[n]^*)\to K_{L,\lambda}(\partial[n-1]^*)$ are given, for $j\in \{0,\ldots,n\}$,  as follows
\begin{equation}\label{pidegboundary}
(\partial_j\psi)_i=  \begin{cases} \psi_i &\mbox{if } i< j \\
\psi_j+\psi_{j+1}  & \mbox{if } i=j\\
\psi_{j+1} & \mbox{if } i>j \end{cases}
\end{equation}
Thus, for instance for $n=3$ one obtains
$$
\left(
\begin{array}{c}
\partial_0 \psi   \\
\partial_1 \psi \\
\partial_2 \psi \\
\partial_3 \psi \\
\end{array}\right)=\left(
\begin{array}{ccc}
 \psi_2 & \psi_3 & \psi_4 \\
 \psi_1+\psi_2 & \psi_3 & \psi_4 \\
 \psi_1 & \psi_2+\psi_3 & \psi_4 \\
 \psi_1 & \psi_2 & \psi_3+\psi_4 \\
\end{array}
\right)
$$
 Hence, for $n>1$   the only solution of the spherical condition $\partial_j\psi=0$ $\forall j$ is $\psi=0$. Indeed, the condition $\partial_0\psi=0$ implies that $\psi_i=0$ for all $i>1$ and the condition $\partial_2\psi=0$ shows that $\psi_1=0$.
\endproof

\subsection{The higher levels of $\bfh(D)$}\label{secthigher}

 One can compute  the  precise number of elements of $\vvert H\Z\vvert_\lambda(k_+)$. We set
\begin{equation}\label{count1}
\gamma(n,k):=\#\left(\vvert H\Z\vvert_n(k_+) \right).
\end{equation}
Any  element $\phi$, with $\sum \vert \phi(x)\vert \leq n$, has a support $s(\phi)\subseteq \{1,\ldots ,k\}$. With the exception of the single element $0$, we let  $m>0$ be the size of the support $s(\phi)=\{x_1,x_2,\ldots, x_m\}$, where the elements are listed in increasing order.  This way, one obtains a list of numbers  all $\leq n$  that is equivalently described by the subset  (with $\psi=\vert\phi\vert$)
$$
\{\psi(x_1),\psi(x_1)+\psi(x_2), \ldots , \psi(x_1)+\psi(x_1)+\ldots + \psi(x_{m})\} \subset \{1,\ldots , n\}.
$$
This process gives, for each of the ${k\choose m}$ subsets with $m$ elements, the number $2^m{n\choose m}$ of possible choices, where the factor $2^m$ takes into account the signs of the $\phi(x)$. Thus one obtains the hypergeometric function
\begin{equation}\label{count2}
\gamma(n,k)=1+\sum_{m=1}^k 2^m{k\choose m}{n\choose m}=\, _2F_1(-k,-n;1;2).
\end{equation}
  Using the vanishing of binomial coefficients, one can replace the index $k$ in  the above sum  by the $\text{inf}(k, n)$  showing the symmetry
\begin{equation}\label{symsym}
\gamma(n,k)=\gamma(k,n).
\end{equation}
Next proposition displays the properties of the higher levels $\pi_n(\bfh(D)(k_+))$

\begin{prop}\label{prophigher} Let  $D=\sum a_j\{p_j\}+ a\{\infty\}$ be a divisor on $\spzb$. The following facts hold
\begin{enumerate}
\item[(i)]~$\pi_1(\bfh(D))=\Vert HL\Vert_{e^a}$,  as $\sss$-modules, with $L=H^0(\Spec \Z,\cO(D_\fii))$.
\item[(ii)]
$
	\#\pi_1(\bfh(D)(k_+)) =\gamma(n,k)$ for $n=\text{IntegerPart}(e^{\deg(D)})$.
	\item[(iii)]~$\pi_0(\bfh(D)(k_+))=* \iff k\leq 2e^{\deg(D)}$.
\item[(iv)]~$\pi_n(\bfh(D)(k_+))$ is trivial for $n>1$.
\end{enumerate}
\end{prop}
\proof $(i)$~The   simplicial set $\bfh(D)(k_+)$ is described in degree $0$, with the notation \eqref{identxi} for the  $\xi_j\in \Hom_\Delta([n],[1])$, $0\leq j\leq n+1$, and with $\lambda =e^a$, by
$$
K_{L,\lambda}(\partial[0]^*\wedge k_+)=K_{L,\lambda}((*,\xi_1)\wedge k_+,(*,\xi_1)\wedge k_+)=(\R/L)^k.
$$
 In degree $1$  it is described by
\begin{equation}\label{pideg1k}
K_{L,\lambda}(\partial[1]^*\wedge k_+)=K_{L,\lambda}((*,\xi_1,\xi_2)\wedge k_+,(*,\xi_2)\wedge k_+)=\Vert H\R \Vert_\lambda(k_+) \times (\R/L)^k.
\end{equation}
The corresponding boundaries  $\partial_j: K_{L,\lambda}(\partial[1]^*\wedge k_+)\to K_{L,\lambda}(\partial[0]^*\wedge k_+)$  are
\begin{equation}\label{pi0high}
\partial_0(\psi)=\psi_2, \ \ \partial_1(\psi)=\psi_1+\psi_2 \qqq \psi=(\psi_1,\psi_2)\in K_{L,\lambda}(\partial[1]^*\wedge k_+).
\end{equation}
The spherical conditions $\partial_j(\psi)=0$ mean that $\psi_2=0\in (\R/L)^k$ and that the image of $\psi_1\in \Vert H\R \Vert_\lambda(k_+)$ is zero in $(\R/L)^k$. Thus  $\psi_1\in \Vert HL \Vert_\lambda(k_+)$. Moreover the relation \eqref{Maydefn} between the spherical elements is trivial for the same reason as stated in the proof of Proposition \ref{proppi0} $(iii)$. \newline
$(ii)$~follows from $(i)$. \newline
$(iii)$~The relation $\cR$ of \eqref{Maydefn} on $K_{L,\lambda}(\partial[0]^*\wedge k_+)=(\R/L)^k$ is given by
$$
(\phi,\phi')\in \cR\iff \exists \psi \in \Vert H\R \Vert_\lambda(k_+)\ \text{s. t.}\ \phi'=\phi +\psi.
$$
Let $d$ be the metric on $\R/L$ which is induced by the usual distance on $\R$. Then
$$
(\phi,\phi')\in \cR\iff \sum d(\phi_j,\phi'_j)\leq \lambda.
$$
Thus  the relation $\cR$ is equal to the square of $(\R/L)^k$, \ie all elements are equivalent   if and only if the sum of the diameters of the metric spaces $\R/L$ is less than $\lambda$, \ie $k/2\leq e^{\deg(D)}$ (one can assume that $D=a\{\infty\}$ so that $L=\Z$ and $\lambda=e^{\deg(D)}$).
\newline
$(iv)$~For $n>1$ one has
\begin{equation}\label{pidegnk}
K_{L,\lambda}(\partial[n]^*\wedge k_+)=\{(\psi_j)_{j\in \{1,\ldots,n\}}\in (\R^k)^{n-1}\times (\R/L)^k\mid \sum_{j\in \{1,\ldots,n-1\}} \Vert \psi_j\Vert \leq \lambda\},
\end{equation}
where the norm of vectors in $\R^k$ is the $\ell^1$-norm.
The boundaries $\partial_j:K_{L,\lambda}(\partial[n]^*\wedge k_+)\to K_{L,\lambda}(\partial[n-1]^*\wedge k_+)$ are given for $j\in \{0,\ldots,n\}$ by
\begin{equation}\label{pidegboundaryk}
(\partial_j\psi)_i=  \begin{cases} \psi_i &\mbox{if } i< j \\
\psi_j+\psi_{j+1}  & \mbox{if } i=j\\
\psi_{j+1} & \mbox{if } i>j \end{cases}
\end{equation}
Thus the same argument of the proof of Proposition \ref{proppi0}	 shows that, for $n>1$ there is no non-zero solution of the spherical condition. \endproof
Proposition \ref{prophigher} $(i)$ shows that the $\sss$-module $\pi_1(\bfh(D))$ qualifies as the module $H^0(\overline{\Spec\Z},\cO(D))$ of global sections of $\cO(D)$. Before  identifying $H^1(\overline{\Spec\Z},\cO(D))$  with  $\pi_0(\bfh(D))$ one needs to take care of the relation $\cR$ defining  the homotopy among  elements of $\pi_0$, since $\cR$ is not in general an equivalence relation. To handle this difficulty one can  implement the
 topos $\frak{Sets}^{(2)}$ introduced in \cite{CCgromov} (and known as the  topos of reflexive (directed) graphs) in order to keep the full  information supplied by  the relation.


\end{document}